\documentclass[12pt,leqno]{article}
\newcommand{\kmcomment}[1]{}
\newcommand{\kmqed}{\hfill\ensuremath{\blacksquare}}

\newcommand{\rank}{\operatorname{rank}}

\newcommand{\RED}[1]{{\color{red}#1}}

\newcommand{\futo}[1]{\mathbf{#1}}

\newcommand{\pdel}{\partial}

\newcommand{\mN}{\ensuremath{\mathbb{N}}} 
\newcommand{\mR}{\ensuremath{\mathbb{R}}} 
\newcommand{\mT}{\ensuremath{\mathbb{T}}} 
\newcommand{\mZ}{\ensuremath{\mathbb{Z}}} 

\newcommand{\frakg}{\mathfrak{g}}

\newcommand{\frakham}{\mathfrak{ham}}
\newcommand{\fraksp}{\mathfrak{sp}}

\newcommand{\HGF}[2]{\text{H}^{#1}_{\rm GF}(\frakham_2^0,\fraksp(2,\mR))_{#2}}

\newcommand{\HGFr}[4]{\text{H}^{#1}_{\rm GF}(\frakham_{#2}^{#3}, \fraksp(#2,\mR))_{#4}}

\newcommand{\Sbt}[2]{[#1,#2]}                
                
\newcommand{\Pkt}[2]{\{#1,#2\}}%
\newcommand{\C}[2]{ C_{#1}^{#2}}       
\newcommand{\CS}[2]{ \text{C}_{#2,#1}} 

\newcommand{\eps}[1]{\epsilon_{#1}}
\newcommand{\mynabla}[2]{\nabla^{#1}_{#2}}

\renewcommand{\RED}[1]{{\color{black}#1}}

\usepackage{theorem}

{
\theoremheaderfont{\bf}
\theorembodyfont{\rmfamily}
\theoremstyle{break}  
\theoremstyle{plain}  

\newtheorem{defn}{\bf Definition}
\newtheorem{prop}{Proposition}[section]
\newtheorem{exam}{Example}[section]
\newtheorem{remark}{Remark}[section]

\theoremstyle{plain}  
\newtheorem{definition}{Definition}[section]
\newtheorem{theorem}{Theorem}[section]

\newtheorem{kmCor}[theorem]{Corollary}
\newtheorem{kmRemark}{Remark}[section]

}

\renewcommand{\[}{$$} \renewcommand{\]}{$$}

\usepackage{setspace,xcolor}
\usepackage[a4paper,truedimen,scale={.85,.85},centering]{geometry}
\usepackage{amsmath,amsfonts,amssymb}
\usepackage[toc,title,page]{appendix}
\usepackage{cases,array}

\usepackage{newtxtext,newtxmath}
\usepackage{fancyvrb,listings,jlisting}
\usepackage{fancyvrb}
\usepackage{amscd}
\newcommand{\san}[2]{f^{[#1]}_{#2}}
\newcommand{\DM}[1]{\text{DM}(#1)}
\newcommand{\sgn}{\operatorname{sgn}}
\newcommand{\myUP}[2]{{#1}^{#2}}

\newcommand{\newPkt}[2]{\Pkt{#1}{#2}_{top}} 
\newcommand{\myZ}[1]{Z_{\ds #1}} 

\newcommand{\myM}[2]{Z\begin{array}{l} #2 \\ #1 \end{array}} 
\newcommand{\ds}{\displaystyle }
\newcommand{\dash}[1]{{#1}'}
\newcommand{\crk}[2]{\operatorname{Corank[1]}(#1,#2)} 

\newcommand{\myMM}[2]{Z[#1]} 

\numberwithin{equation}{section}

\title{
First Betti number of 
weighted homology group of 
Hamiltonian vector fields on symplectic tori}

\author{
Hiroki Kodama\thanks{ 
Graduate School of Mathematical Sciences,
University of Tokyo, 3-8-1 Komaba,
Meguro-ku, Tokyo 153-8914, JAPAN
kodama@ms.u-tokyo.ac.jp
}\and 
Kentaro Mikami\thanks{
  Department of Computer Science and Engineering
    Akita University, partially supported by Grant-in-Aid for
      Scientific Research (C) of JSPS, JP26400063 and JP23540067}
       \and Tadayoshi Mizutani\thanks{
       Professor Emeritus, Saitama University }
       }

\date{March 2018 (Second version)}

\begin{document}\parindent=0pt

\allowdisplaybreaks
\maketitle
\thispagestyle{empty}

\section{Introduction} 
Given a Lie algebra \(\ds \frakg\), we have  \RED{the
chain complex 
$(\ds \Lambda^m\frakg, \pdel)$,}  \RED{ with 
the boundary operator on the decomposable elements  given by} 
\[\ds \pdel (\xi_{1}\wedge \cdots \wedge
\xi_{m}) = \sum_{i<j} (-1)^{i+j} [\xi_{i},\xi_{j}] \wedge \xi_{1} 
\wedge \cdots \wedge  \widehat{\xi_{i}} 
\wedge \cdots \wedge \widehat{\xi_{j}}  
\wedge \cdots \wedge \xi_{m}\] 
where \(\ds \xi_{1},\ldots, \xi_{m} \in \frakg\) and 
\(\ds \widehat{\xi_{i}}\) means omitting  
\(\ds \xi_{i}\).  
Then the \RED{$m$-th homology group of the complex is} defined by  
\[\ds \ker ( \pdel : \Lambda ^m
\frakg \longrightarrow \Lambda ^{m-1}\frakg )/ \operatorname{image} ( \pdel :
\Lambda ^{m+1} \frakg \longrightarrow \Lambda ^{m}\frakg ) \;,\] 
\RED{and this is the homology group of $\frakg$ with coefficient $\mR$ 
with trivial action.}
\RED{Dually, one can define cohomology groups by considering 
the linear functionals on $\frakg$}.  
If \(\ds \frakg\) is infinite dimensional, then we collect
continuous linear functionals as the dual space, with respect to some
topology, and we obtain a continuous cohomology groups of \(\ds \frakg\).  

As introduced in \cite{KM:affirm}, two sort of 
well-studied examples of infinite dimensional Lie algebras are Lie
algebras of volume preserving vector fields on \(\ds \mR^{n}\) and 
(formal) Hamiltonian vector fields on the symplectic \(\ds \mR^{2n}\). 
Those cohomology  groups are sometimes called Gel'fand-Fuks
cohomologies.  

There is an interesting geometric observation by 
Kontsevich (\cite{Kont:RW}).  
Let $\ds {\cal F}$ be a foliation on a manifold $M$. We have the
foliated cohomology defined by $\ds \text{H}_{\cal F}^{\bullet}
(M,\mR) := \text{H}^{\bullet}( \Omega_{\cal F})$ where $\ds
\Omega_{\cal F} = \Omega(M)/\text{I}({\cal F})$, $\ds\Omega(M)$ is
the exterior algebra of differential forms on $M$, and
$\ds\text{I}({\cal F})$ is the ideal generated by $\ds\{ \sigma
\in \Omega^{1}(M) \mid \langle \sigma, \text{T}{\cal F}\rangle =0
\}$.   
Kontsevich (\cite{Kont:RW}) showed that 
if $\ds {\cal F}$ is a codimension $2n$ foliation endowed with a symplectic
form $\omega$ in the direction
 transverse to  the foliation,  then there is a commutative
diagram: 
\begin{equation}\begin{CD}
        @. \text{H}_{\cal F}(M,\mR) @>\omega^{n}\wedge>> \text{H}_{\text{
        DR}}^{\bullet+2n}(M,\mR)  @.  \\
@.   @AAA   @AAA \\
@.   \HGFr{\bullet}{2n}{0}{ } @>\omega^{n}\wedge>> \HGFr{\bullet+2n}{2n}{ }{ } @.
\end{CD}
\label{Kont:CD}
\end{equation} 
where   $\ds\frakham_{2n}^{0}$ is the Lie subalgebra of the
Hamiltonian vector fields of the formal polynomial vanishing at
the origin of $\ds  \mR^{2n}$.  
Metoki (\cite{metoki:shinya}) 
 found the next non-trivial group 
$\ds  \HGF{9}{14}{} $ next to Gel'fand-Kalinin-Fuks's work. 
Kotschick and Morita (\cite{KOT:MORITA}) determined the
space $\ds  \HGF{\bullet}{w}{} $  for $w\le 10$. 
Mikami-Nakae-Kodama  
(\cite{M:N:K}) we determined  $\ds\HGF{\bullet}{w}{0} $
for $w\le 22$. 

Concerning Kontsevich homomorphism given in
the bottom line of (\ref{Kont:CD}), 
there are works by 
Kotschick and Morita (\cite{KOT:MORITA}) and Mikami
(\cite{KM:affirm}). 

In many area of geometry, research target next to Euclidean spaces seems
to be tori. Thus, in this paper we develop 
the Lie algebra (co)homology theory of Hamiltonian vector fields 
on the symplectic torus $\ds \mT^n = \mR^n/\mZ^n$ ($n$ even).  
Main result we got is: 

\textbf{Theorem:}\quad 
Let \(\ds \crk{w}{2n}\) be the dimension of the first homology group of \(\ds
 \CS{w}{1} \mathop{\leftarrow}^{\pdel } \CS{w}{2}\), where $w$ is the
 weight. Then we have the next recursive  formula:  
\begin{align*} 
 \crk{w}{2n} 
=&  4 + 4 \sum_{i=1}^{w}  \crk{i}{2n-2} - 3 \crk{w}{2n-2} \;, 
\\
\noalign{and}
\crk{0}{2n} =& 1 \; .
\end{align*}

\section{Recall of symplectic $\mR ^{n}$ case ($n$ is even)}
With respect to the Poisson bracket of 
the standard 
symplectic $\ds \mR^{n}$, the polynomial algebra has a Lie algebra
structure. So, we may study its Lie algebra 
(co)homology groups.  However, the (co)spaces are huge, we introduce the
notion of weight to reduce the (co)chain spaces. 
For each multi-index 
\(\ds A \in \mN^{n}\) of length $n$, we use the notation 
\(\ds 
 z[A] \) for the monomial \(\ds  
 x_1^{A[1]} \cdots x_{n}^{A[n]}\).   
The degree of \(z[A]\) is defined by
\(\ds |A| := \sum_{i=1}^{n} A[i]\). 

Consider the vector space $\ds V_{h}$ spanned by $\ds\{ z[A] \mid |A|
= h \}$, namely the subspace of $h$-homogeneous polynomials. 
The dimension of \(\ds V_{h}\) is known as $\tbinom{n-1+h}{n-1}$.

We define the weight of 
$h$-homogeneous polynomial by \( h-2 \). A reason is:   
For $f$-homogeneous polynomial $F$ and 
$h$-homogeneous polynomial $H$, take  the Poisson bracket $\Pkt{F}{H}$. 
The degree of  $\Pkt{F}{H}$ is $f+h-2 = (f-2) + (h-2) +2 $, and   
we have $wt(\Pkt{F}{H}) = wt(F)+wt(H) $ for weights. 

In Lie algebra homology theory, $m$-th chain space is given by 
\[ z[A^{(1)}] \wedge z[A^{(2)}] \wedge \cdots \wedge z[A^{(m)}] \quad
\text{with}\quad |A^{(j)}|>0\quad\text{for}\quad j=1..m\] 
and the 
$0$-th chain space is spanned by \( z[A]\) with \( |A|=0\) and $\mR$
itself. 

We define the weight of 
\(\ds z[A^{(1)}] \wedge z[A^{(2)}] \wedge \cdots \wedge z[A^{(m)}]\)
naturally by 
\[
 wt(z[A^{(1)}] \wedge z[A^{(2)}] \wedge \cdots \wedge z[A^{(m)}]) :=
 wt(z[A^{(1)}])+ wt(z[A^{(2)}]) +\cdots +  wt(z[A^{(m)}]) \]
 and we may consider the subspace  \(\ds \CS{w}{m} \) of 
$m$-th chain space with the weight $w$. 

 \(\ds \CS{w}{m} \) 
is spanned by
\begin{align} & \sum \Lambda^{k_1} V_1 \otimes  \Lambda^{k_2} V_2 \otimes 
    \cdots \otimes 
 \Lambda^{k_{\ell}} V_{\ell} \label{base:1} \\ 
\noalign{where} 
& k_1+k_2+ \cdots + k_{\ell} = m\;, \label{base:2}\\  
 & (1-2)k_1 + (2-2)k_2+ (3-2)k_3 + \cdots + (\ell-2)k_{\ell} = w\;,
 \label{base:3}\\ 
 \noalign{and}   
& 0 \leq k_j \leq \tbinom{n-1+j}{n-1} \;.  
\label{base:4}
\end{align} 
Using \eqref{base:2} and \eqref{base:3}, we have 
\begin{align}
& 1k_1+ 2k_2 + \cdots + \ell k_{\ell} = w + 2m \; . \label{base:5}
\end{align}
As observed in \cite{M:N:K}, \eqref{base:2} and \eqref{base:5} means 
the sequence $[k_1,\ldots,k_{\ell}]$ corresponds to a  Young
diagram of area $w+2m$ with length $m$ as follows:  

Starting \eqref{base:5} with $k_{\ell}\ne 0$, the sequence $\ds ( 
\underbrace{\ell,\ldots,\ell}_{k_{\ell}}, \ldots, 
\underbrace{1,\ldots,1}_{k_{1}})$ is a Young diagram of length m and the
area is $ w+2m$.  

Conversely, take a Young diagram $\lambda$ with area $w+2m$ and length m.   
Let $\ds \mu$ be the conjugate Young diagram of $\lambda$.  In short, 
$\ds \mu_{i} = \text{the depth of } i\text{-th column of }\lambda$.  
Let $\ell := \text{length}(\mu) = \lambda_{1}$ and 
put $\ds \mu_{1+\ell} = 0$, and   
define \[k_i := \mu_{i} - \mu_{1+i} \quad\text{for}\quad  i=1..\ell\;.\]  
Then 
\begin{align*}
\sum k_i &= (\mu_1-\mu_2)+(\mu_2-\mu_3)+\cdots +(\mu_{\ell} -0) =
\mu_{1} = \text{length}(\lambda) \\
\sum i k_i &= (\mu_1-\mu_2)+2 (\mu_2-\mu_3)+\cdots + \ell (\mu_{\ell} -0) =
\mu_{1} + \mu_{2} + \cdots + \mu_{\ell} = \text{area}(\mu) = 
\text{area}(\lambda) \end{align*} 
Subtracting 2-times of \eqref{base:2} from \eqref{base:5}, we have that 
$ -k_1 + k_3 + 2k_4 + \cdots = w$ and so we see that $w \geq -k_1 \geq -n$ using
\eqref{base:4}. 
Subtracting \eqref{base:2} from \eqref{base:5}, we see that 
$w+2m \geq m $, i.e., $ m \geq \max(-w,0)$.  
Subtracting 3-times of \eqref{base:2} from  \eqref{base:5}, we have that 
$ -2 k_1 - k_2 + k_4 + 2k_5 + \cdots = w-m$ and so we see that 
$w -m +2k_1+k_2  \geq 0$, i.e.,  
$m \leq w +2k_1+k_2 \leq w + n(n+5)/2$  
using \eqref{base:4}. We summarize these relations:
\[ w \geq -n\quad\text{and}\quad  \max(-w,0) \leq m \leq w + n(n+5)/2\;. \]

\begin{remark}
Consider a constant Poisson tensor field 
\(\ds \pi  = \sum_{i=1}^{k} \frac{\pdel}{\pdel x_{2i-1} }
\wedge \frac{\pdel}{\pdel x_{2i} }\) 
\kmcomment{
\(\ds \pi  = \sum_{i=1}^{k} {\pdel}/ {\pdel x_{2i-1} }
\wedge \pdel/{\pdel x_{2i} }\) 
}
on \(\ds \mR^n\) with \(2k\leq n\). 
Now the homogeneity is 0 and 
the Euler characteristic of Lie algebra homology groups of the Poisson
bracket on \(\ds\mR^n\) is 0 for each weight $w$ (cf.\ 
\cite{Mik:Miz:homogPoisson}). 

\end{remark}

\kmcomment{
\subsection{Young diagrams (general theory)}
In the previous subsection, we knew that the $m$-th (co)chain space with
the weight $w$ corresponds with the set of Young diagrams of length $m$
and area $w+2m$, we use the notation \(\ds \mynabla{w+2m}{m}\). Here we
neglect dimensional restriction, and use our notation of Young diagrams.   
\(k = [k_{1},\ldots,k_{\ell}] \in \mynabla{A}{m} \) if and only if 
\(\ds \sum k_{j} = m \) and 
\(\ds \sum j k_{j} = A \).  
We would like to add one cell to a Young diagram 
\(k = [k_{1},\ldots,k_{\ell}] \in \mynabla{A}{m} \).  

One way is to put the cell on the most left bottom of \(k\), i.e., 
if we say the new Young diagram by \(\ds\dash{k}\), then 
\begin{equation}  \dash{k}_{1} = 1+ k_{1}\quad \text{and}\quad    
\dash{k}_{j} = k_{j}\quad  \text{for}\quad  j>1\;.\quad\text{Now}\quad   
\dash{k} \in \mynabla{A+1}{m+1}\;. \end{equation} 

The other way is keeping the length and attach the
one cell to the right end of some row. If this is possible, say the
width of the row by $j$. Then $\ds k_{j}$ is reduced by 1 and     
 $\ds k_{j+1}$ increases by 1, namely      
\begin{equation} \dash{k}_{j} = -1+ k_{j}\quad\text{and}\quad    
\dash{k}_{j+1} = 1+ k_{j+1}\;.\quad\text{The possibility is}\quad  k_{j}
>0 \quad \text{ for some}\quad j\;. \label{eqn:add:one} \end{equation} 
The reverse operation,  deleting one cell from a given Young diagram, is
also characterized by our notation of Young diagrams as below.  
\begin{align}
\ell =  [\ell_{1},\ell_{2},\ldots ] \in \mynabla{A}{m}\quad  
& 
\text{if}\quad  \ell_{1} > 0 \quad\text{then}\quad 
 [ -1+\ell_{1},\ell_{2},\ldots ] \in \mynabla{A-1}{m-1} \\
& \text{if}\quad \ell_{j} > 0\; (\exists j>1)\quad\text{then}\quad   
[ \ldots, 1+\ell_{j-1},-1+\ell_{j}, \ldots ] \in \mynabla{A-1}{m}  
\end{align}
For \(k \in \mynabla{a}{m}\), we have a subset, say 
\(\ds \Hat{k}\) 
of \(\ds \mynabla{a+1}{m}\) by the operation \eqref{eqn:add:one}. 
Take \(\ell \in \mynabla{a}{m}\)
and we have 
\(\Hat{\ell}\in \mynabla{a+1}{m}\). 
Suppose $k \ne \ell$. What is 
\(\ds \Hat{k}\cap  \Hat{\ell}\) ?

If \(\ds k_{i} \geq 1\), \(\ds \ell_{i+1}\geq 2\), 
\(\ds k_{i+2}\geq 1\),
\(\ds k_{i} = \ell_{i}+1\), 
\(\ds k_{i+1} = \ell_{i+1}-2\), 
\(\ds k_{i+2} = \ell_{i+2}+1\) 
and \(\ds k_{j}= \ell_{j}\) otherwise,  then 
\[\begin{array}{*{4}{c}l} 
[\ldots, &  -1+k_{i}, & 1+k_{i+1}, & k_{i+2},&  \ldots ] \\
= [\ldots , & \ell_{i}, &  -1 +\ell_{i+1}, & 1+ \ell_{i+2}, & \ldots ]
\in \Hat{k}\cap\Hat{\ell}\;.\end{array}  
\]

For some $i < j-1$,  
if \(\ds k_{i} \geq 1\), \(\ds \ell_{j}\geq 1\), 
\(\ds k_{i} = \ell_{i}+1\), 
\(\ds k_{i+1} = \ell_{i+1}-1\), 
\(\ds k_{j} = \ell_{j}-1\),  
\(\ds k_{j+1} = 1+\ell_{j}\),  
and \(\ds k_{s}= \ell_{s}\) otherwise,  then 
\[ \begin{array}{*{6}{c}l}  
[\ldots, &  -1+k_{i}, & 1+k_{i+1}, & \ldots, & k_{j}, & k_{j+1}, &\ldots ]
\\
=[\ldots,& \ell_{i}, & \ell_{i+1}, &\ldots, & -1+ \ell_{j}, &
1+\ell_{j+1}, &
\ldots ]
\in \Hat{k}\cap\Hat{\ell}\;.\end{array}  
\]

\begin{exam}

\[
\begin{array}{lc@{\quad\rightarrow\quad}l}
\mynabla{4}{2} = \{ [1,0,1], [0,2] \} 
& 
$k$ &  \Hat{k} \\
& 
$[1,0,1]$ & $[0,1,1], [1,0,0,1]$ \\
& 
$[0,2]$ & $[0,1,1]$ \\\hline
\mynabla{5}{2} = \{ [1,0,0,1], [0,1,1,0] \} 
& 
$k$ &  \Hat{k} \\
& 
$[1,0,0,1]$ & $[0,1,0,1], [1,0,0,0,1]$ \\
& 
$[0,1,1]$ & $[0,0,2], [0,1,0,1]$ 
\end{array} 
\]

\end{exam}
}

\section{Preparation on $\ds \mT^n$}

\subsection{Functions on $\ds\mT^{n}$}
Functions on 
$\ds \mT^n = \mR^n/(2\pi Z)^n $ are considered as 
periodic functions on $\ds \mR^n$ of the period $2\pi$.  
Instead of generators 
$\ds x_1,\ldots, x_n$ of polynomial algebra of \(\ds \mR^n\), 
we may deal with the algebra generated by 
\begin{equation}  \cos x_1,\;  \sin x_1,\;  \ldots,\;  \cos x_n,\;  \sin x_n
\label{eqn:base}
\end{equation} 

\begin{defn}
For a given non-negative integer $k$, 
define the space 
$W_{k}$ spanned by \begin{equation} \{ 
\sin^{a_1} x_1 \sin^{a_2} x_2  \cdots \sin^{a_n} x_n  
\cos^{a_1} x_1 \cos^{a_2} x_2  \cdots \cos^{a_n} x_n  \mid 
\sum_{i=1}^{n}(a_i+b_i) \leq k\}\;, \label{eqn:famous:formula}\end{equation} 
\end{defn}
Then, we have a filtration in the algebra generated by \eqref{eqn:base}
as follows:
\[ W_{0} \subset W_{1}  \subset W_{2}  \subset \cdots \]

\begin{remark}
Even if 
we require the equality "$=$" in
\eqref{eqn:famous:formula} in the definition above, 
by the well-known formula  $\ds \cos ^2 \theta + \sin^2 \theta=1$,  
$W_{2}$ contains constant function, and
$W_{1} \subset W_{3}$, $W_{0} \cup W_{2} \subset 
W_{4}$, and so on.  \end{remark}

\subsection{Concrete basis} \label{subsec:concrete}
\begin{prop}
We use the notation 
\begin{equation} z[a_1,\ldots,a_n,b_1,\ldots, b_n] :=  
\sin^{a_1} x_1 \cdots \sin^{a_n} x_n \cdot 
\cos^{b_1} x_1 \cdots \cos^{b_n} x_n \label{eqn:nota:a}\;. \end{equation} 

For $k\in\mN$, 
\(\ds 
z[a_1,\ldots,a_n,b_1,\ldots, b_n] \) with \(\ds\sum_{i=1}^{n} (a_i + b_i)
=k\),
\(\ds a_j \in \mN\) 
and \(\ds b_j \in\{0,1\}\ \text{for}\ j=1,\ldots, n\) form a basis of the 
quotient space 
$\ds W_{k}/W_{k-1} $ and    
\begin{equation} \dim 
(W_{k}/W_{k-1} ) 
= \sum_{\alpha+\beta=k} \tbinom{n}{\alpha} \tbinom{n-1+\beta}{n-1}
\;.
\label{nota:a:dim}
\end{equation}

Consider the subspace of $\ds W_{k}$ spanned by the above basis and
denote the subspace by $\ds V_{k} $.  Then 
$\ds V_{k} \cong W_{k}/W_{k-1} $ and we see that  $\ds W_{\infty} =
\mathop{\oplus}_{k=0}^{\infty} V_{k}$.  
\end{prop}

\begin{remark}
In the academic meeting at Tambara, Japan in fall 2017 we presented this
topic, and got a suggestion by Professor Tsuboi 
that it is better to use another basis coming from 
\begin{equation} 
\sin {a_1} x_1 \cos {b_1} x_1 
\cdots \sin {a_n} x_n \cos {b_n} x_n \label{eqn:nota:b} \end{equation} 
instead of \eqref{eqn:nota:a}.  
When we started studying this topic we are faced  which basis to choose. 
At that time, we chose  \eqref{eqn:nota:a} because that is more similar
to ``polynomial algebra''.  In the appendix, we expand our discussion using
another basis \eqref{eqn:nota:b}.  

\end{remark}

\subsection{Multiplication on $\ds V_{i} \times V_{j}$}
Take 
$\ds z[A,B] \in V_{i}$ and $\ds z[A',B'] \in V_{j}$. 
\begin{align*} z[A,B] z[A',B'] =& z[P,Q]\; \text{where}\ P=A+A', Q = B + B' 
\; (\text{naturally first, then }) \\
    =& (-1)^{ k_1+k_2+\cdots+k_n } z[P',Q'] 
\end{align*}
with \(\ds     
    P'[\ell] = P[\ell]+2 k_{\ell}\), \(\ds 0\leq Q'[\ell]=
    Q[\ell]-2k_{\ell}< 2\) for $\ell$.   

\subsection{Naturally reduced Poisson bracket}
In this subsection, we do not assume $n$ is even. We may consider 
a Poisson tensor given by 
\[ \pi = \sum p_{ij} \frac{\pdel}{\pdel x_i} \wedge  \frac{\pdel}{\pdel
x_j} \] where   $\ds p_{ij}\in V_{h}$, and  
\(\ds \Pkt{\cdot}{\cdot}\) is the bracket defined by $\pi$. 
For \(\ds F \in V_{f}\) and  \(\ds J \in V_{j}\),  
it holds 
\( \Pkt{F}{J} \in W_{f+j+h} \).  
Making use of decomposition into $\ds V_{k}$, 
we have 
\[\ds\Pkt{F}{J} = \Pkt{F}{J}_{top} + \Pkt{F}{J}_{lower} \in
V_{f+j+h}\oplus W_{f+j+h-1}\;. \]
Since 
\[\ds\Pkt{ \Pkt{F}{J}}{H}_{top} = \Pkt{\Pkt{F}{J}_{top}}{K}_{top} \]
holds and so 
\[ (F,J) \mapsto \Pkt{F}{J}_{top} \] satisfies the Jacobi identity.  
\(\ds \Pkt{V_{f}}{V_{j}}_{top} \subset V_{f+j+h}\) holds from the
definition. 

\begin{definition}
For each 
\( F \in V_{f} \),  define the its weight by $\ds wt(F):= f+h$, where $h$
is the homogeneous degree of the Poisson tensor. 
\end{definition}
Then the next holds: \[ wt(\Pkt{F}{J}_{top} ) = wt(F) + wt(J) \quad 
\text{for}\ F\in V_{f}\;,\;  J\in V_{j}\;.  
\] 
Thus, we use $\ds \Pkt{F}{J}_{top}$ as a new bracket: 
$\ds\newPkt{F}{J} := \text{the top term of} \Pkt{F}{J}$. Then it induces
a Lie algebra structure on 
$\ds \frakg = \mathop{\oplus}_{j=0}^{\infty} V_j $. 
We shall study 
Lie algebra (co)homology groups with respect to
this new bracket. 
 $m$-th chain space is given by  
\[\C{m}{} = \sum \Lambda^{k_1}V_{1} \otimes \cdots \otimes 
\Lambda^{k_\ell}V_{\ell} \]
with 
$\ds k_1 + \cdots + k_{\ell} = m$.  

Take care of weight, 
\[ \CS{w}{m} = \sum \Lambda^{k_1}V_{1} \otimes \cdots \otimes 
\Lambda^{k_\ell}V_{\ell} \]
with 
$\ds k_1 + \cdots + k_{\ell} = m$ and   
$\ds k_1 + 2 k_2 + \cdots + \ell k_{\ell} = w - h m $, 
and correspond to Young diagrams of area 
$w - h m$, height $m$, and there is a dimensional restriction  
\(\ds 0 \leqq k_j \leqq \sum_{a+b=j} \tbinom{n-1+a}{n-1} \tbinom{n}{b} \).

For \(\ds b_j \in\{0,1\}\ \text{for}\ j=1..n\),    
\begin{align*}
& 
\frac{\pdel}{\pdel x_j} 
z[a_1,\ldots,a_n,b_1,\ldots, b_n]\\
=& 
a_j z[a_1,\ldots,a_j- 1,\ldots,b_1,\ldots, b_j+1,\ldots ]  
- b_j z[a_1,\ldots,a_j + 1,\ldots,b_1,\ldots, b_j-1,\ldots ]  
\\
=& \begin{cases}
a_j z[a_1,\ldots,a_j- 1,\ldots,b_1,\ldots, b_j+1,\ldots ]  
 & \text{if } b_j=0 \\
- (a_j+b_j) z[a_1,\ldots,a_j + 1,\ldots,b_1,\ldots, b_j-1,\ldots ]  
& \text{if } b_j = 1
\end{cases}
\\
=& (-1)^{b_j} (a_j + b_j) 
 z[a_1,\ldots,a_j - (-1)^{b_j} ,\ldots,b_1,\ldots, b_j + (-1)^{b_j} ,\ldots ]  
 \\
 =& 
(-1)^{b_j} (a_j + b_j) 
 z[a_1,\ldots,a_j +2b_j -1 ,\ldots,b_1,\ldots, 1-b_j ,\ldots ]  
\end{align*}
is a formula of basic differentiation, we may calculate the Poisson
bracket by this formula and multiplication rule. 

When we use multi-index set of length $2n$, the next holds 
\begin{align} 
\frac{\pdel}{\pdel x_j} 
z[C ] 
& = (-1)^{C_{n+j} } (C_{j} + C_{n+j})\; z[ C' ]  \label{so:1}\\ 
\noalign{where}C_{n+j} \in \{0,1\}\quad\text{and}\quad  
C'_{j} &= C_{j} + 2C_{n+j} -1\; ,\quad 
C'_{n+j} = 1- C_{n+j} \;. \label{so:2}\;.  \end{align}

\section{2-dim symplectic torus}
In this section, $n=2$ and consider the standard Poisson tensor and the
standard Poisson bracket on 2-dimension torus.  

\subsection{Basic result} \label{subsec:main}

\begin{theorem}\label{thm:first}
On the symplectic 2-torus \(\ds \mT^2\), we consider the formal
Hamiltonian vector fields, we may identify with periodic ``polynomials''. 
We saturate the space with respect to the weight and have $m$-th chain
space $\ds \CS{w}{m}$ with the given weight $w$.  
We have a chain complex 
\[  0 \leftarrow \CS{w}{1} 
 \mathop{\leftarrow}^{\pdel} \CS{w}{2} 
 \mathop{\leftarrow}^{\pdel} \CS{w}{3}  
 \mathop{\leftarrow}^{\pdel} \cdots \]
The first Betti number of this chain complex is 4 for each positive
weight $w$. 
\end{theorem}

\textbf{Proof:} 
Since $\ds 0 \leftarrow \CS{w}{1} $ holds  
for even $n$-dimensional standard symplectic torus, 
in order to know the first Betti number, we need data of 
dimension of $\ds\CS{w}{1}$  and 
$\ds\CS{w}{2}$, $a,b$, and  
the kernel dimension of 
$\ds 
\CS{w}{1} \mathop{\leftarrow}^{\pdel} \CS{w}{2} $ 
$c$ (or the rank), then it turns out $\text{1st Betti} = a-b+c $ 
as the next table shows: 
Actually, since  
\(\ds \text{coker}(\pdel) := \CS{w}{1}/ \pdel (\CS{w}{2}) \), 
\(\ds \text{corank}(\pdel) = \dim 
 \text{coker}(\pdel) = \dim \CS{w}{1} - \dim  \pdel (\CS{w}{2}) = 
  \dim \CS{w}{1} - \rank \pdel  =  a+c-b  \). 
\[
\begin{array}{c|*{6}{c}}
m & 0 & \leftarrow & \ds \CS{w}{1} & \leftarrow &  \CS{w}{2}  \\\hline
\dim & && a  &  & b   \\
\ker & && a &  & c  \\\hline
\text{Betti} &&& a+c-b  & 
\end{array}
\]
Here we introduce the notation 
$\crk{w}{n}$ 
the corank at degree 1 
of the boundary operator \(\pdel\) to \(\ds \CS{w}{1} = V_{w,n}\)  
where \(V_{w,n}\) means the $w$-homogeneous "polynomial" space on
\(\mT^{n}\) for even $n$. 

For a given weight
$w$, $\ds \CS{w}{1} = V_{w,n}$ and  
$\ds \CS{w}{2} $ is given as follows depending to the parity of $w$.  

\[ 
\CS{w}{2} = \begin{cases}\ds 
\sum_{j=1}^{(w-1)/2} V_{j,n} \otimes V_{w-j,n} & \text{if}\; 
w\; \text{odd} \\
\ds
\sum_{j=1}^{(w-1)/2} V_{j,n} \otimes V_{w-j,n} \oplus 
\Lambda^2 V_{w/2,n} & \text{if}\; w\; \text{even} \end{cases} \]

Now back to $n=2$. 
Using \eqref{so:1} and  \eqref{so:2}, 
we express concretely the bracket 
$\ds \Pkt{z[X]}{z[Y]}$, where $X$ and $Y$ are length $2n=4$ multi-index sets. 
Here now we change our notation \(\ds z[A]\) by 
\(\ds \myM{a_1,a_2}{a_3,a_4}\) where \(\ds a_3, a_4\) run only 0 or 1.   
With this notation, our Poisson bracket for generators (before
``modulation'') is given by
\begin{align}  \Pkt{\myM{A}{P}}{\myM{B}{Q}}  
=& |A,B| \myM{A+B - (\eps{1}+\eps{2})} {P+Q + (\eps{1}+\eps{2})} 
+ |P,Q| \myM{A+B + (\eps{1}+\eps{2})} {P+Q - (\eps{1}+\eps{2})} 
\\& 
- |A,Q| \myM{A+B + (\eps{1}-\eps{2})} {P+Q - (\eps{1}-\eps{2})} 
- |P,B| \myM{A+B + (-\eps{1}+\eps{2})} {P+Q - (-\eps{1}+\eps{2})} 
\\
\noalign{taking care of ``modulation''} 
=& |A,B| \myM{A+B - (\eps{1}+\eps{2})} {P+Q + \eps{1}+\eps{2}} 
+ |P,Q| \myM{A+B - (\eps{1}+\eps{2})} {P+Q +\eps{1}+\eps{2}} 
\\& 
+ |A,Q| \myM{A+B - (\eps{1}+\eps{2})} {P+Q + \eps{1}+\eps{2}} 
+ |P,B| \myM{A+B - (\eps{1}+\eps{2})} {P+Q + \eps{1}+\eps{2}} 
\\
=& \left( |A,B| + |A,Q| + |P,B| + |P,Q| \right) 
\myM{A+B - [1,1]} {P+Q + [1,1]} \label{mk:eqn}
\end{align} 
where $|A,B|$ means the determinant of 2 by 2 matrix \(\begin{bmatrix}
A\\ B\end{bmatrix}\), $\eps{1}=[1,0]$ and  $\eps{2}=[0,1]$.  

When $P=Q$, 
\eqref{mk:eqn} implies 
\begin{align*}  \Pkt{\myM{A}{P}}{\myM{B}{P}}  
=& \left( |A,B| + |A,P| + |P,B| + |P,P| \right) 
\myM{A+B - [1,1]} {2P + [1,1]} \;. 
\\\noalign{In particular, putting
$P=[0,0]=O$, then we have} 
\Pkt{\myM{A}{O}}{\myM{B}{O}}  
=&  |A,B|  \myM{A+B - [1,1]} {[1,1]} 
=  (a_{1}b_{2}-a_{2}b_{1}) \myM{A+B - [1,1]} {[1,1]} 
\end{align*} 
Putting \(\ds a_{2}=b_{1}=0, a_{1} -1 = c_{1}, b_{2}-1=c_{2}\), we have 
\begin{align*}  \Pkt{\myM{{[1+c_{1},0]}}{O}}{\myM{{[0,1+c_{2}]}}{O}}  
=&  (1+c_{1}) (1+c_{2})  \myM{{[c_{1},c_{2}]}}{[1,1]} 
\end{align*} 
\begin{align} 
 \myM{{[c_{1},c_{2}]}} {[1,1]} 
=&  \frac{1}{(1+c_{1}) (1+c_{2})} 
 \Pkt{\myM{{[1+c_{1},0]}}{O}}{\myM{{[0,1+c_{2}}]}{O}}  
 \label{eqn:11}
\end{align} 
When $P \neq Q$, assume $p_{2} \neq q_{2}$.  
Let \(\ds 
\hat{P}= [p_{1}, q_{2}], 
\hat{Q}= [q_{1}, p_{2}]\).  
Using 
\eqref{mk:eqn}, we have another equation:  
\begin{align}  \Pkt{\myM{A}{\hat{P}}}{\myM{B}{\hat{Q}}}  
=& \left( |A,B| + |A,\hat{Q}| + |\hat{P},B| + |\hat{P},\hat{Q}| \right) 
\myM{A+B - [1,1]} {\hat{P}+\hat{Q} + [1,1]} 
\label{eqn:hat}
\end{align} 
Since \(\ds \hat{P} + \hat{Q} = P+Q =[p_{1}+q_{1},1] \), subtracting
\eqref{eqn:hat} from \eqref{mk:eqn}, we have 
\begin{align*}  
\Pkt{\myM{A}{P}}{\myM{B}{Q}} -  
\Pkt{\myM{A}{\hat{P}}}{\myM{B}{\hat{Q}}}  
=& 
\left( |A,Q| + |P,B| + |P,Q| 
 - |A,\hat{Q}| - |\hat{P},B| - |\hat{P},\hat{Q}| \right) 
\myM{A+B - [1,1]} {[p_{1}+q_{1}+1, 2]} 
\end{align*} 
The coefficient of the right hand side of the above equation is \(\ds
-(p_{2}-q_{2}) (a_{1}+b_{1} + p_{1} + q_{1})\).  
So, we have 
\begin{align*}  
& \Pkt{\myM{{[a_{1},a_{2}]}}{P}}{\myM{{[b_{1},b_{2}]}}{Q}} -  
\Pkt{\myM{{[a_{1},a_{2}]}}{\hat{P}}}{\myM{{[b_{1},b_{2}}]}{\hat{Q}}}  
\\
=& -(p_{2}-q_{2}) (a_{1}+b_{1} + p_{1} + q_{1})  
\myM{{[a_{1}+b_{1}-1, a_{2}+b_{2}-1]}} {{[p_{1}+q_{1}+1, 2]}} 
\\
=& 
(p_{2}-q_{2}) (a_{1}+b_{1} + p_{1} + q_{1})  
\myM{{[a_{1}+b_{1}-1, a_{2}+b_{2}+1]}} {{[p_{1}+q_{1}+1, 0]}} 
\\ 
\noalign{thus,} 
& \Pkt{\myM{{[1+a_{1},a_{2}]}}{P}}{\myM{{[b_{1},b_{2}]}}{Q}} -  
\Pkt{\myM{{[1+a_{1},a_{2}]}}{\hat{P}}}{\myM{{[b_{1},b_{2}]}}{\hat{Q}}}  
\\ 
=& 
(p_{2}-q_{2}) (1+a_{1}+b_{1} + p_{1} + q_{1})  
\myM{{[a_{1}+b_{1}, 1+a_{2}+b_{2}]}}{{[p_{1}+q_{1}+1, 0]}} 
\end{align*} 
Therefore, 
\begin{subequations}
\begin{align} 
& (p_{2}-q_{2}) 
\myM{{[a_{1}+b_{1}, 1+a_{2}+b_{2}]}}{{[p_{1}+q_{1}+1, 0]}} 
\\ =& \frac{1}{1+a_{1}+b_{1} + p_{1} + q_{1} } 
\left( 
\Pkt{\myM{{[1+a_{1},a_{2}]}}{P}}{\myM{{[b_{1},b_{2}]}}{Q}} -  
\Pkt{\myM{{[1+a_{1},a_{2}]}}{\hat{P}}}{\myM{{[b_{1},b_{2}]}}{\hat{Q}}}  
\right)
\end{align} 
\end{subequations}
When \(\ds p_{1}+q_{1}=0\), we have  
\begin{subequations}
\begin{align} 
& (p_{2}-q_{2}) 
\myM{{[a_{1}+b_{1}, 1+a_{2}+b_{2}]}}{{[1, 0]}} \label{eqn:10}\\
=& \frac{1}{1+a_{1}+b_{1} + p_{1} + q_{1} } 
\left( 
\Pkt{\myM{{[1+a_{1},a_{2}]}}{P}}{\myM{{[b_{1},b_{2}]}}{Q}} -  
\Pkt{\myM{{[1+a_{1},a_{2}]}}{\hat{P}}}{\myM{{[b_{1},b_{2}]}}{\hat{Q}}}  
\right)
\end{align} 
\end{subequations} 
When \(\ds p_{1}+q_{1}=1\), we have  
\begin{subequations}
\begin{align} 
& 
(p_{2}-q_{2}) 
\myM{{[a_{1}+b_{1}, 1+a_{2}+b_{2}]}}{{[2, 0]}} 
\notag \\
=& \frac{1}{2+a_{1}+b_{1} } 
\left( 
\Pkt{\myM{{[1+a_{1},a_{2}]}}{P}}{\myM{{[b_{1},b_{2}]}}{Q}} -  
\Pkt{\myM{{[1+a_{1},a_{2}]}}{\hat{P}}}{\myM{{[b_{1},b_{2}]}}{\hat{Q}}}  
\right) \notag\\
\noalign{and so, after modulation, 
    replacing \(\ds a_{1}+1\) by \(a_{1}\), we have } 
& (p_{2}-q_{2}) 
\myM{{[1+a_{1}+b_{1}, 1+a_{2}+b_{2}]}}{{[0, 0]}} 
\label{eqn:00:by10}
\\
=& \frac{-1}{1+a_{1}+b_{1} } 
\left( 
\Pkt{\myM{{[a_{1},a_{2}]}}{P}}{\myM{{[b_{1},b_{2}]}}{Q}} -  
\Pkt{\myM{{[a_{1},a_{2}]}}{\hat{P}}}{\myM{{[b_{1},b_{2}]}}{\hat{Q}}}  
\right)
\end{align} 
\end{subequations} 
When $P \neq Q$ and $p_{1} \neq q_{1}$, we follow the same discussion
above and get   
\begin{subequations}
\begin{align}
& (p_{1}-q_{1}) 
\myM{{[1+a_{1}+b_{1}, a_{2}+b_{2}]}} {{[0,p_{2}+q_{2} + 1]}} 
\\
=& \frac{-1}{ 1+a_{2}+b_{2} + p_{2} + q_{2} } 
\left( 
\Pkt{\myM{{[a_{1},1+a_{2}]}}{P}}{\myM{{[b_{1},b_{2}]}}{Q}} -  
\Pkt{\myM{{[a_{1},1+a_{2}]}}{\tilde{P}}}{\myM{{[b_{1},b_{2}]}}{\tilde{Q}}}  
\right)
\end{align} 
\end{subequations}
where 
\(\ds \tilde{P}= [q_{1}, p_{2}]\) and \(\ds  
\tilde{Q}= [p_{1}, q_{2}]\).  

\begin{subequations}
\begin{align}
\noalign{ If \(\ds p_{2}+q_{2}=0\)}  
& (p_{1}-q_{1}) 
\myM{{[1+a_{1}+b_{1}, a_{2}+b_{2}]}} {{[0, 1]}} 
\label{eqn:01}
\\ =& \frac{-1}{ 1+a_{2}+b_{2} } 
\left( 
\Pkt{\myM{{[a_{1},1+a_{2}]}}{P}}{\myM{{[b_{1},b_{2}]}}{Q}} -  
\Pkt{\myM{{[a_{1},1+a_{2}]}}{\tilde{P}}}{\myM{{[b_{1},b_{2}]}}{\tilde{Q}}}  
\right)
\end{align} 
\end{subequations}

\begin{subequations}
\begin{align}
\noalign{ If \(\ds p_{2}+q_{2}=1\), then }  
& (p_{1}-q_{1}) 
\myM{{[1+a_{1}+b_{1}, a_{2}+b_{2}]}} {{[0, 2]}} 
\notag \\
=& \frac{-1}{ 2+a_{2}+b_{2}  } 
\left( 
\Pkt{\myM{{[a_{1},1+a_{2}]}}{P}}{\myM{{[b_{1},b_{2}]}}{Q}} -  
\Pkt{\myM{{[a_{1},1+a_{2}]}}{\tilde{P}}}{\myM{{[b_{1},b_{2}]}}{\tilde{Q}}}  
\right)
\notag\\
\noalign{and so, we have}
& (p_{1}-q_{1}) 
\myM{{[1+a_{1}+b_{1}, 1+a_{2}+b_{2}]}} {{[0, 0]}} 
\label{eqn:00:by01}
\\
=& \frac{1}{ 1+a_{2}+b_{2}  } 
\left( 
\Pkt{\myM{{[a_{1},a_{2}]}}{P}}{\myM{{[b_{1},b_{2}]}}{Q}} -  
\Pkt{\myM{{[a_{1},a_{2}]}}{\tilde{P}}}{\myM{{[b_{1},b_{2}]}}{\tilde{Q}}}  
\right) 
\end{align} 
\end{subequations}

We summarize the above discussion in a table:

\begin{center}
\begin{tabular} {*{3}{|c}|}
\hline
type in \(\ds V_{w}\)  & bracket image &  cokernel of bracket \\ \hline
\(\ds \myM{c_1,c_2}{0,0} \) & 
    c.f.\ \eqref{eqn:00:by01} or \eqref{eqn:00:by10}  & 
\(\ds \myM{w,0}{0,0}\) and \(\ds \myM{0,w}{0,0}\)
\\\hline
\(\ds \myM{c_1,c_2}{0,1} \) & c.f.\ \eqref{eqn:01} & 
\(\ds \myM{0,w-1}{0,1}\) 
\\\hline
\(\ds \myM{c_1,c_2}{1,0}\) & c.f.\ \eqref{eqn:10} & 
\(\ds \myM{w-1,0}{1,0}\) 
\\\hline
\(\ds \myM{c_1,c_2}{1,1} \) & c.f.\ \eqref{eqn:11} & none\\
\hline
\end{tabular}
\end{center}
This concludes the corank of  \(\ds \CS{w}{1}
\mathop{\leftarrow}^{\pdel} \CS{w}{2} \) is 4, namely the first Betti
number is 4 for each positive weight $w$.  
\kmqed

\begin{kmRemark}
For the standard symplectic space \(\ds \mR^2 \), the Poisson bracket is
given by 
\[ \Pkt{\futo{x}^{A}}{\futo{x}^{B}} = |A B| \futo{x}^{A+B-[1,1]} \]
where \(\ds \futo{x}^{A} = x_{1}^{a_{1}}  x_{2}^{a_{2}}\).  
Putting \(\ds a_{2}=b_{1}=0\), then replacing 
\(\ds a_{1}-1\) by \(\ds a_{1}\) and 
\(\ds b_{2}-1\) by \(\ds b_{2}\), we have  
\[ \Pkt{\futo{x}^{a_{1}+1,0}}{\futo{x}^{0,b_{2}+1}} = 
(1+a_{1}) (1+b_{2}) 
\futo{x}^{a_{1}, b_{2} } \]
This says the first Betti number is always zero for every weight in the
case of standard symplectic \(\ds \mR^{2}\).  We may prove that for general even
dimensional case the first Betti number is zero by apply our discussion
in this article.   
\end{kmRemark}

\subsection{Concrete examples of $n=2$ case}
The next several tables are list of the dimensions of chain spaces and
the kernel spaces and Betti numbers for the lower weights. 

\begin{tabular}{c|*{3}{r}}
wt=2 & 1 & $\leftarrow$ 2 \\\hline
dim & 8 & 6 \\
ker & 8 & 2 \\\hline
Betti & 4 &2
\end{tabular}
\hfil 
\begin{tabular}{c|*{3}{r}}
wt=3 & 1 & $\leftarrow$ 2 & $\leftarrow$ 3 
\\\hline
dim & 12 & 32 & 4 \\
ker & 12 & 24 & 0 \\\hline
Betti & 4 & 20 & 0
\end{tabular} 
\hfil 
\begin{tabular}{c|*{4}{r}}
wt=4 & 1 & $\leftarrow$ 2 & $\leftarrow$ 3 & $\leftarrow$ 4 
\\\hline
dim & 16 & 76 & 48 & 1 \\
ker & 16 & 64 & 10 & 0 \\\hline
Betti & 4 & 26 & 9 & 0
\end{tabular}

\begin{tabular}{c|*{4}{r}}
wt=5 & 1 & $\leftarrow$ 2 & $\leftarrow$ 3 & $\leftarrow$ 4 
\\\hline
dim & 20 & 160 & 184 & 32 \\
ker & 20 & 144 & 76 & 0 \\\hline
Betti & 4 & 36 & 44 & 0
\end{tabular}
\hfil  
\begin{tabular}{c|*{5}{r}}
wt=6 & 1 & $\leftarrow$ 2 & $\leftarrow$ 3 & $\leftarrow$ 4 & $\leftarrow$ 5 
\\\hline
dim & 24 & 274 & 536 & 216 & 8 \\
ker & 24 & 254 & 324 & 16 & 0 \\\hline
Betti & 4 & 42 & 124 & 8 & 0
\end{tabular} 

\section{The first Betti numbers on general symplectic tori}
\subsection{The first Betti numbers when $n=4$}

\begin{tabular}{c|*{3}{r}}
wt=2 & 1 & $\leftarrow$ 2 \\\hline
dim & 32 & 28 \\
ker & 32 & 20 \\\hline
Betti & 24 & 20
\end{tabular}
\hfil 
\begin{tabular}{c|*{3}{r}}
wt=3 & 1 & $\leftarrow$ 2 & $\leftarrow$ 3 
\\\hline
dim & 88 & 256 & 56 \\
ker & 88 & 208 & 16 \\\hline
Betti & 40 & 168 & 16
\end{tabular} 
\hfil 
\begin{tabular}{c|*{4}{r}}
wt=4 & 1 & $\leftarrow$ 2 & $\leftarrow$ 3 & $\leftarrow$ 4 
\\\hline
dim & 192 & 1200 & 896 & 70 \\
ker & 192 & 1064 & 436 & 4 \\\hline
Betti & 56 & 604 & 330 & 4
\end{tabular}

\begin{tabular}{c|*{5}{r}}
wt=5 & 1 & $\leftarrow$ 2 & $\leftarrow$ 3 & $\leftarrow$ 4  & $\leftarrow$ 5 
\\\hline
dim & 360 & 4352 & 6432 & 1792 & 56 \\
ker & 360 & 4064 & 3816 & 304 & 0 \\\hline
Betti & 72 & 1448 & 2328 & 248 & 0
\end{tabular}
\hfil  
\begin{tabular}{c|*{6}{r}}
wt=6 & 1 & $\leftarrow$ 2 & $\leftarrow$ 3 & $\leftarrow$ 4 & $\leftarrow$ 5  & $\leftarrow$ 6 
\\\hline
dim & 608 & 12852 & 32864 & 18816 & 2240 & 28 \\
ker & 608 & 12332 & 23304 & 5488 & 80 & 0 \\\hline
Betti & 88 & 2772 & 9976 & 3328 & 52 & 0
\end{tabular}

\medskip

Pick up only the first Betti numbers. Then we have 
\begin{tabular}{c|*{5}{c}}
wt & 2 & 3 & 4 & 5 & 6 \\\hline
1st Betti & 24 & 40 & 56 & 72 & 88 \\\hline
difference &   & 16 & 16 & 16 & 16
\end{tabular} and 
so we {\color{black}expect} the first Betti numbers consists of 
arithmetic progression with common difference 16. 
In deed, 
we have the next theorem.  
\begin{theorem}\label{thm:dim4}
The second difference of the sequence of 
the first Betti number 
of the weight $w$ on $\ds \mT^{4}$ 
is an arithmetic progression with common difference 16.  
More precisely, 
the first Betti number for the weight $w>0$ is given by  
\( 4^2 w -8 \).
\end{theorem}
\textbf{Proof:} 
Denote 
\(\ds  z_{\ds u_1,u_2,u_3,u_4, u_5,u_6,u_7,u_8} \) by 
\(\ds 
\myM{{[u_1,u_2,u_3,u_4]}}{{[u_5,u_6,u_7,u_8]}}\), and  
\(\ds \myM{A,A'}{P,P'}\) in short.  
The Poisson bracket \(\ds \Pkt{\myM{A,A'}{P,P'}}{\myM{B,B'}{Q,Q'}} \) is 
given as 
\begin{align*}
& \Sbt{ \Sbt{
\frac{\pdel}{\pdel x_1} \wedge \frac{\pdel}{\pdel x_2} + 
\frac{\pdel}{\pdel x_3} \wedge \frac{\pdel}{\pdel x_4} }{
 \myM{A,A'}{P,P'}}}  {\myM{B,B'}{Q,Q'}} \\
=&  
\Sbt{ \Sbt{
\frac{\pdel}{\pdel x_1} \wedge \frac{\pdel}{\pdel x_2}}{
  \myM{A,A'}{P,P'}}}  {\myM{B,B'}{Q,Q'}} 
+   
\Sbt{ \Sbt{
\frac{\pdel}{\pdel x_3} \wedge \frac{\pdel}{\pdel x_4} }{
  \myM{A,A'}{P,P'}}}  {\myM{B,B'}{Q,Q'}} \\ 
=&  
\Sbt{ \Sbt{
\frac{\pdel}{\pdel x_1} \wedge \frac{\pdel}{\pdel x_2}}{
  \myM{A}{P}}}  {\myM{B}{Q}} \myM{A'+B'}{P'+Q'}
+   
 \myM{A+B}{P+Q}
\Sbt{ \Sbt{
\frac{\pdel}{\pdel x_3} \wedge \frac{\pdel}{\pdel x_4} }{
  \myM{A'}{P'}}}  {\myM{B'}{Q'}}  
 \end{align*} 
Thus we have \kmcomment{using the notation \(\ds [1,1] = \Delta\),} 
\begin{align}
 \Pkt{\myM{A,A'}{P,P'}}{\myM{B,B'}{Q,Q'}} \label{dim4:mk:master}
\equiv & 
 \Pkt{\myM{A}{P}}{\myM{B}{Q}}^{(2)} \myM{A'+B'}{P'+Q'}
 + \myM{A+B}{P+Q} \Pkt{\myM{A'}{P'}}{\myM{B'}{Q'}}^{(2)} 
\end{align} 
where \(\ds {\Pkt{\cdot}{\cdot}}^{(2)}\) means the Poisson bracket in
2-dimensional torus.  
In \eqref{dim4:mk:master}, put \(\ds B'=Q'=O\) or \(\ds B=Q=O\).   
Then because of 
 \(\ds \myM{B'}{Q'} \) or \(\ds \myM{B}{Q} \)  is constant, we have 
\begin{align}
 \Pkt{\myM{A,A'}{P,P'}}{\myM{B,O}{Q,O}} 
\equiv & \Pkt{\myM{A}{P}}{\myM{B}{Q}}^{(2)} \myM{A'}{P'} \label{cok:one}\\
\Pkt{\myM{A,A'}{P,P'}}{\myM{O,B'}{O,Q'}} 
\equiv & 
 \myM{A}{P} \Pkt{\myM{A'}{P'}}{\myM{B'}{Q'}}^{(2)} \label{cok:two}
\end{align} 

\kmcomment{
Here we introduce the notation 
$\crk{w}{n}$ 
the corank 
of the boundary operator \(\pdel\) to \(V_{w,n}\)
where \(V_{w,n}\) means the $w$-homogeneous "polynomial" space on
\(\mT^{n}\) for even $n$. 
}

We know that $\crk{0}{2}=1$ and $\crk{w}{2}=4$ for $w>0$.   
\eqref{cok:one} and \eqref{cok:two} imply the diagram
\[
\begin{array}{c||c|c} 
V_{i,2} \backslash V_{j,2} & \text{out of image in }\mT^2 & \text{image in
}\mT^2 \\\hline\hline
\text{out of image in }\mT^2 & \text{out of image in }\mT^4 
& \text{image in }\mT^4 \\\hline
\text{image in }\mT^2 & \text{image in }\mT^4 & \text{image in } \mT^4
\end{array}
\]
and we see that for $w>0$ 
\begin{align*} 
 \crk{w}{4} 
=& 
 \crk{0}{2} \crk{w}{2} \\& 
 + \sum_{i=1}^{w-1} \crk{i}{2} \crk{w-i}{2} 
 + 
  \crk{w}{2} \crk{0}{2} 
  \\
=&  4 + 4^2 (w-1) + 4 = 4^2 w -8 
\\
\crk{0}{4} =& 1 
\end{align*}
\kmqed

\subsection{The first Betti numbers for general even $n$}

Observing the discussion in \(\ds \mT^4\) case, we show the same
discussion works well for general even dimensional torus \(\ds
\mT^{2n}\) case, too.  
\begin{theorem} \label{thm:gen}
Let \(\ds \crk{w}{2n}\) be the dimension of the first homology group of \(\ds
 V_{w,2n} \mathop{\leftarrow}^{\pdel } \CS{w}{2}\). Then we have the
 next recursive  formula:  
\begin{align*} 
 \crk{w}{2n} 
=&  4 + 4 \sum_{i=1}^{w}  \crk{i}{2n-2} - 3 \crk{w}{2n-2} 
\\
\crk{0}{2n} =& 1 
\end{align*} 
\end{theorem}
\textbf{Proof:} 
Let \(\ds \myM{\myUP{A}{1}\ldots \myUP{A}{n}}
              {\myUP{P}{1}\ldots \myUP{P}{n}}\) be our basis.  
To avoid complicated expression, we abbreviate the first $(n-1)$
sequences 
\(\ds \myUP{A}{1}\ldots \myUP{A}{n-1}\) by \(A\) and 
\(\ds \myUP{P}{1}\ldots \myUP{P}{n-1}\) by \(P\).   
Then the
              Poisson bracket satisfies 
\begin{subequations}
\begin{align} 
\Pkt{ 
 \myM{A\;\myUP{A}{n}}
              {P\;\myUP{P}{n}}} 
{ \myM{B\; \myUP{B}{n}}
              {Q\;\myUP{Q}{n}}}  
=& 
\Pkt{ 
 \myM{A} {P}}
{ \myM{B }{Q} } 
 \myM{\myUP{A}{n}+\myUP{B}{n} }
              {\myUP{P}{n}+ \myUP{Q}{n}} 
\label{shiki:1}
              + 
 \myM{A+B} {P +Q } 
 \Pkt{ \myM{\myUP{A}{n}}{\myUP{P}{n}}}
{ \myM{\myUP{B}{n}}{\myUP{Q}{n}}}
 \end{align}
\end{subequations}
\(\ds \myUP{B}{n}= \myUP{Q}{n}= O\) implies 
\begin{subequations}
\begin{align} 
\Pkt{ 
 \myM{A\;  \myUP{A}{n}}
              {P\; \myUP{P}{n}}} 
{ \myM{B\; O} {Q\; O}} 
=& 
\Pkt{ 
 \myM{A } {P }} 
{ \myM{B } {Q }} 
 \myM{\myUP{A}{n} }
              {\myUP{P}{n} } 
\label{shiki:2:2}
\end{align}
\end{subequations}
Now replace \(\ds \myUP{P}{n}\) by 
\(\ds \myUP{P}{n} + \myUP{Q}{n}\), 
and \(\ds \myUP{A}{n}\) by 
\(\ds \myUP{A}{n} + \myUP{B}{n}\), 
we have 
\begin{subequations}
\begin{align} 
&
\Pkt{ \myM{A\; \myUP{A}{n}+\myUP{B}{n} } 
          {P\; \myUP{P}{n}+\myUP{Q}{n} }} 
    { \myM{B \;  O} {Q\; O}}  
= 
\Pkt{ 
  \myM{A } {P }} 
{ \myM{B } {Q }}
 \myM{\myUP{A}{n}+\myUP{B}{n} 
 }
              {\myUP{P}{n}+\myUP{Q}{n} 
              } 
\label{shiki:3}
 \end{align}
\end{subequations}
Subtracting \eqref{shiki:3} from \eqref{shiki:1}, we have 
\begin{subequations}
\begin{align}
&
\Pkt{ 
 \myM{A\; \myUP{A}{n}}
              {P\;  \myUP{P}{n}}} 
{ \myM{B\; \myUP{B}{n}}
              {Q\; \myUP{Q}{n}}} -  
\Pkt{ \myM{A\;  \myUP{A}{n}+\myUP{B}{n} }
              {P \; \myUP{P}{n}+\myUP{Q}{n} }} 
{ \myM{B\; O} {Q\; O}} 
= 
 \myM{A+B } {P+Q } 
\Pkt{ \myM{\myUP{A}{n}}{\myUP{P}{n}}}
{ \myM{\myUP{B}{n}}{\myUP{Q}{n}}} \label{shiki:4:2} 
 \end{align}
\end{subequations}

\eqref{shiki:2:2} and  
\eqref{shiki:4:2} imply the diagram
\[
\begin{array}{c||c|c} 
V_{i,2n-2} \backslash V_{j,2} & \text{out of image in }\mT^2 & \text{image in
}\mT^2 \\\hline\hline
\text{out of image in }\mT^{2n-2} &  
\text{out of image in }\mT^{2n} &  
 \text{image in }\mT^{2n} \\\hline
\text{image in }\mT^{2n-2} & \text{image in }\mT^{2n} & \text{image in }
\mT^{2n}
\end{array}
\]
and we see that for $w>0$ 
\begin{align*} 
 \crk{w}{2n} 
=& 
 \crk{0}{2n-2} \crk{w}{2} \\& 
 + \sum_{i=1}^{w-1} \crk{i}{2n-2} \crk{w-i}{2} 
 \\& 
 + \crk{w}{2n-2} \crk{0}{2} 
  \\
=&  4 + 4 \sum_{i=1}^{w-1}  \crk{i}{2n-2} + \crk{w}{2n-2} 
\\
\crk{0}{2n} =& 1 
\end{align*}
\kmqed

\begin{kmCor} 
\label{cor:dim6}
On \(\ds \mT^6\), 
 \(\ds \crk{0}{6}  = 1\) and for \(\ds w>0\) we have   
\begin{equation}  
 \crk{w}{6}  = 4 ( 8 w^2 - 12 w + 7)\; .
\label{eqn:dim6} 
 \end{equation} 

On \(\ds \mT^{2n+2}\), it holds 
\begin{equation}\ds 
\sum_{k=0}^ n (-1)^{k} \binom{n}{k} 
 \crk{w-k}{2n+2}  = 2^{2n+2} \label{eqn:recursive} \end{equation}  
for \( w > n\). 
\end{kmCor}
\textbf{Proof:} 
Combining Theorem \ref{thm:dim4} and Theorem \ref{thm:gen}, direct
computation yields the formula 
\eqref{eqn:dim6}. Also,  
\begin{align*}
& 
\sum_{k=0}^{2} (-1)^{k} \binom{2}{k} 
 \crk{w-k}{6}  \\=&  
 \crk{w}{6}  - 2   
 \crk{w-1}{6}  +   
 \crk{w-2}{6}   = 2^6\;.  
\end{align*}
Thus, \eqref{eqn:recursive} holds for $n=2$.  
Assume \eqref{eqn:recursive} holds on \(\ds \mT^{2n+2}\) for general $n$.   
First we confirm that 
\begin{align} & 
\sum_{k=0}^{n+1} (-1)^k \binom{n+1}{k} \crk{w-k}{2n+2}
\label{eqn:hojo}
\\ =& 
\sum_{k=0}^{n+1} (-1)^k ( \binom{n}{k} 
+\binom{n}{k-1})  
\crk{w-k}{2n+2} = 2^{2n+2}-2^{2n+2} = 0\;.  
\notag
\end{align}
Now we prove 
\eqref{eqn:recursive} holds on \(\ds \mT^{2n+4}\).  
\begin{align*} & 
\sum_{k=0}^{n+1} (-1)^k \binom{n+1}{k} \crk{w-k}{2n+4} \\
=& \sum_{k=0}^{n+1} (-1)^k \binom{n+1}{k} \left( 4+ 4 \sum_{i=1}^ {w-k}
\crk{i}{2n+2} - 3 \crk{w-k}{2n+2} \right) 
\\ 
=& 4 \sum_{k=0}^{n+1} (-1)^k \binom{n+1}{k}  - 3 
 \sum_{k=0}^{n+1} (-1)^k \binom{n+1}{k} 
\crk{w-k}{2n+2} \\& 
+ 
 4 \sum_{k=0}^{n+1} (-1)^k \binom{n+1}{k}  \sum_{i=1}^ {w-k}
\crk{i}{2n+2} \\
=& 
 4 \sum_{k=0}^{n+1} (-1)^k \binom{n+1}{k}  \sum_{i=1}^ {w-k}
\crk{i}{2n+2} \\
\noalign{because of general fact about alternating sum, and 
\eqref{eqn:hojo}. Changing the order of summation, we have }
=& 4 \sum_{i=1}^{w-n-1} \sum_{k=0}^{n+1} (-1)^k
\binom{n+1}{k}\crk{i}{2n+2} 
+ 4 \sum_{i=w-n}^{w} \sum_{k=0}^{w-i} (-1)^k \binom{n+1}{k}
\crk{i}{2n+2}\\
=& 0  
+ 4 \sum_{i=w-n}^{w} \sum_{k=0}^{w-i} (-1)^k \binom{n+1}{k}
\crk{i}{2n+2}\\
=&   4 \sum_{i=w-n}^{w} \sum_{k=0}^{w-i} (-1)^k ( \binom{n}{k} +
 \binom{n}{k-1} ) \crk{i}{2n+2}
 \\
=&   4 \sum_{i=w-n}^{w} (-1)^{w-i}  \binom{n}{w-i} \crk{i}{2n+2}
 = 4\; 2^{2n+2}
 = 2^{2n+4} \;. 
\end{align*}
\kmqed

\kmcomment{
\subsubsection{Examples of the first Betti numbers when $n=6$}

Again, we just take the 
the first Betti numbers for lower weight. 
\begin{center}
\begin{tabular}{c|*{5}{c}}
wt & 2 & 3 & 4 & 5 & 6 \\\hline
1st Betti & 60 & 172 & 348 & 588 & 892  \\\hline
difference & - & 112 & 176 & 240 & 304 \\\hline
difference  & - & - & 64  & 64 & 64 \\\hline
\end{tabular}
\end{center}
}

\subsection{Not symplectic but Poisson structures on tori}
In 
\cite{Mik:Miz:homogPoisson} we develop Lie algebra (co)homology theory
for homogeneous Poisson structures on $\mR^n$.  
On $\mR^n$ we also consider   
homogeneous Poisson structures using our basis. A typical 0-homogeneous
Poisson structure is \(\ds 
\frac{\pdel}{\pdel x_{1}} \wedge \frac{\pdel}{\pdel x_{2}} + 
\cdots + 
\frac{\pdel}{\pdel x_{2m-1}} \wedge \frac{\pdel}{\pdel x_{2m}}\) on
\(\ds \mT^{2m + (n-2m)}\).  We denote our basis as 
\( \myM{A,A'}{P,P'} \) where  
\( \myM{A}{P} \) are basis of \(\ds \mT^{2m }\) and 
\( \myM{A'}{P'} \) are basis of \(\ds \mT^{n- 2m }\).  
Then the Poisson bracket is given by 
\[ \Pkt{ \myM{A,A'}{P,P'} } { \myM{B,B'}{Q,Q'} } = 
\Pkt{ \myM{A}{P} } { \myM{B}{Q} } 
\myM{A'+B'}{P'+Q'}\;. \] 
Like Theorem \ref{thm:gen}, 
on the first Betti number of the above Poisson structure, 
we may have a formula 
\begin{align*} 
 \crk{w}{n} 
=&  \sum_{i=0}^{w}  \crk{i}{2m} \dim V_{w-i, n-2m} 
\;, 
\end{align*} 
where  \(\ds \dim V_{k, n-2m} = 
\sum_{\alpha+\beta=k} \tbinom{n-2m-1+\alpha}{n-2m - 1}
\tbinom{n-2m}{\beta}
 \).
 
\bibliographystyle{plain}
\bibliography{km_refs}

\begin{appendices}

\section{New basis}
On $n$-torus, 
\(\ds W_{k}\) is generated by at most $k$ product of 
$\ds \cos x_1, \sin x_1, \ldots, \cos x_n, \sin x_n$. We apply
Fourier expansion, we have another basis for
\(\ds W_{k} / W_{k-1} \) defined by the followings: 
\begin{align} \myZ{a_1,\ldots,a_n} =  
& \san{1}{a_1}  \cdots 
 \san{n}{a_n} \quad  
 \text{where}\quad \sum_{i=1}^{n} |a_i| = k \quad\text{and}\quad 
\san{j}{a}  = 
 \begin{cases} \sin a x_j & \text{if}\quad a > 0 \\
 \cos a x_j & \text{if}\quad a \leq 0 \end{cases} \label{ap:nota:b}
\end{align}
\subsection{Dimension formulae}
\begin{prop} \label{prop:new:basis:dim}
Let \(P_{k}\) be the vector spaced linearly spanned by \eqref{ap:nota:b}. 
Then 
\begin{equation} 
\dim P_{k} = \sum_{\ell} \tbinom{n}{\ell}\tbinom{k-1}{ \ell-1} 
2^{\ell} \quad\text{for}\quad  k>0\;,\quad\text{and}\quad \dim P_{0} = 1
\label{eqn:new} 
\end{equation} 
\end{prop}
Proof of Proposition \ref{prop:new:basis:dim}: 
\kmcomment{
$\hrulefill$

整数列 \(\ds a_{1},\ldots, a_{n}\) 但し
\(\ds \sum_{i} | a_{i} | = k\) 
を数え上げます。
$n$個の整数の列
\(\ds a_{1},\ldots, a_{n}\) 
に対し, $\ds a_{i} \ne 0 $ なるものの個数を $\ell$ とする。その場所は
$\ds \tbinom{n}{\ell}$ 通りであり, 長さ $\ell$ の正整数(絶対値を付けた)
の和が $k$ となるのは, 1 を均等配分した後, $k-\ell$ 個の玉を$\ell$個の皿
に分配する仕方なので, \(\ds \tbinom{\ell-1+k-\ell}{\ell-1}\) 通りある。
以上から \eqref{eqn:new} を得る。

$\hrulefill$
}
We count the number of sequences of integers 
\(\ds a_{1},\ldots, a_{n}\) of length $n$ 
with \(\ds \sum_{i} | a_{i} | = k\).  
For a given 
\(\ds a_{1},\ldots, a_{n}\),  
let $\ell$ be the number of non-zero 
$\ds a_{i}$'s.  Those patterns are of 
$\ds \tbinom{n}{\ell}$ cases. 
The possibility of length $\ell$ sequences of {\color{red}positive} 
integers with 
$k$ total sum  is known as 
\(\ds \tbinom{\ell-1+k-\ell}{\ell-1}\).  
Thus we have \eqref{eqn:new}.  
\kmqed


\newcommand{\Dimm}[2]{\text{Dim}_{#1}(#2)}
\begin{prop}
Although  
\eqref{nota:a:dim} and \eqref{eqn:new} are the cardinalities of bases
which look different, 
the two values are equal because  
the space is the same.  Namely, 
\[
 \sum_{\alpha+\beta=k} \tbinom{n}{\alpha} \tbinom{n-1+\beta}{n-1}
 = \sum_{\ell} \tbinom{n}{\ell}\tbinom{k-1}{ \ell-1} 
2^{\ell} 
 \] 
\end{prop}
\textbf{Proof:} 
We put the dimension of $P_{w}$ on $n$-torus by \( \ds \DM{n,w}\).  
For $n>1$  
\begin{equation}
\DM{n,w} =  \sum_{i+j=w} \DM{1,i}  \DM{n-1,j} 
\label{key:eqn}
\end{equation}
holds, and 
\( \ds \DM{1,w} = \begin{cases} 1 & w=0 \\
2 & w > 0 \end{cases}\quad \)   
implies 
\begin{equation}
\label{key2:eqn}
\DM{n,w}   
  = \DM{n-1,w} + 2 \sum_{j=1}^{w} \DM{n-1, j-1}\;. 
  \end{equation}
Reducing $w$ by 1 in \eqref{key2:eqn},  we have 
\[
\DM{n,w-1}   
  = \DM{n-1,w-1} + 2 \sum_{j=1}^{w-1} \DM{n-1, j-1} \;.
\] 
Subtracting the above from \eqref{key2:eqn}, we have  
\begin{equation}
\label{key3:eqn}
\DM{n,w} - \DM{n,w-1}  
  = \DM{n-1,w} + \DM{n-1, w-1}\;. 
  \end{equation}
Under the condition \(\DM{n,0} = 1\),  
\eqref{key3:eqn} implies
\eqref{key2:eqn}, in other words, 
\eqref{key2:eqn} and \eqref{key3:eqn} are equivalent. 

Denoting 
\eqref{nota:a:dim} by \(\ds \Dimm{1}{n,k}\), and  \eqref{eqn:new} by  
\(\ds \Dimm{2}{n,k}\),   
we shall show 
\(\Dimm{1}{n,w}\) and \(\Dimm{2}{n,w}\) satisfy the same 
relation  \eqref{key3:eqn},  
we complete the proof of 
\(\ds \Dimm{1}{n,w} = \Dimm{2}{n,w} = \DM{n,w} \).   
Hereafter, we use the formula 
\(\ds \tbinom{n}{\ell} =  \tbinom{n-1}{\ell} +
\tbinom{n-1}{\ell-1} \) frequently.

\paragraph{Old basis:}
\begin{align*}
 \Dimm{1}{n,w} 
=& \sum \tbinom{n}{\ell} \tbinom{n-1+w-\ell}{n-1} = 
\sum \tbinom{n}{\ell} (\tbinom{n-2+w-\ell}{n-2} +  \tbinom{n-2+w-\ell}{n-1}) 
\\
=& 
\sum (  \tbinom{n-1}{\ell} +\tbinom{n-1}{\ell-1} ) 
\tbinom{n-2+w-\ell}{n-2} +  
\sum \tbinom{n}{\ell} \tbinom{n-2+w-\ell}{n-1} 
\\ = &
\Dimm{1}{n-1,w} + 
\Dimm{1}{n-1,w-1}   + 
\sum \tbinom{n}{\ell} \tbinom{n-2+w-\ell}{n-1} 
\\
{\color{red}
=} & 
\Dimm{1}{n-1,w} + 
\Dimm{1}{n-1,w-1} + 
\Dimm{1}{n,w-1} \;. 
\end{align*}
Thus, we have the same relation 
\begin{equation} 
\Dimm{1}{n,w} - \Dimm{1}{n,w-1} =  
\Dimm{1}{n-1,w} + \Dimm{1}{n-1,w-1} 
\end{equation}
as \eqref{key:eqn}.   
\kmcomment{
Sum up the equations above, we have the same relation with 
\eqref{key2:eqn}.  
}
\kmqed

\paragraph{New Basis:}
\begin{align}
\Dimm{2}{n,w} =&  \sum \tbinom{n}{\ell} \tbinom{w-1}{\ell-1} 2^{\ell} 
=    
\sum (\tbinom{n-1}{\ell} + \tbinom{n-1}{\ell-1} ) 
\tbinom{w-1}{\ell-1} 2^{\ell} 
\notag \\
=&  \Dimm{2}{n-1,w} + 2 \sum \tbinom{n-1}{\ell-1}  
\tbinom{w-1}{\ell-1} 2^{\ell-1} \label{AA} \\
= & 
\Dimm{2}{n-1,w} + 2 \sum \tbinom{n-1}{\ell-1}  
( \tbinom{w-2}{\ell-2} + 
\tbinom{w-2}{\ell-1}) 2^{\ell-1} \notag  \\
= & 
\Dimm{2}{n-1,w} + 2 \Dimm{2}{n-1,w-1 } + 
2 \sum \tbinom{n-1}{\ell-1}  
\tbinom{w-2}{\ell-1} 2^{\ell-1} \notag \\
{\color{red} =} &  
\Dimm{2}{n-1,w} + 2 \Dimm{2}{n-1,w-1 }
 + \left( \Dimm{2}{n,w-1} - \Dimm{2}{n-1,w-1}\right)  \label{BB}
\end{align}
where we used the relation \eqref{AA} to get \eqref{BB}. 
Again, we have the same relation as \eqref{key:eqn}.  
\kmqed

\textbf{Alternative Proof:}
The space \(\ds W_{k}/W_{k-1} \) explained in 
\eqref{nota:a:dim} of Subsection \ref{subsec:concrete}
has dimension 
\( \ds  \sum_{\alpha+\beta=k} \tbinom{n}{\alpha} \tbinom{n-1+\beta}{n-1}
\)  and this number is the cardinality of
\[\ds D_{1}(n,k) = \{ 
(a_1,\ldots,a_n,b_1,\ldots, b_n)\in \mN^{n} \times \{0,1\}^{n}  \mid \sum_{i=1}^{n} (a_i + b_i) = k\; \}\;. \]  
On the other hand, the space 
\(P_{k}\) spanned by 
the new basis comes from 
\[\ds D_{2}(n,k) = \{ 
(c_1,\ldots,c_n)\in \mZ^{n} \mid \sum_{i=1}^{n} |c_i | = k\; \}\;. \]  
We define a map $\psi$ by 
\[ \psi( a_{i}, b_{i} ) = \begin{cases}
a_{i} + 1 & \text{if}\; b_{i} = 1\;, \\
-a_{i}  & \text{if}\; b_{i} = 0 \;.\end{cases} \] 
Using $\psi$, we define another map
\[\ds \Psi : D_{1}(n,k) \ni 
(a_1,\ldots,a_n,b_1,\ldots, b_n)\mapsto 
(\psi(a_1,b_1), \ldots, \psi(a_n, b_n)) \in
D_{2}(n,k)\;.\] 
\(\ds \Psi\)  is well-defined and bijective.  Thus, 
\[\ds 
\Dimm{1}{n,k} = \# D_{1}(n,k) = \# D_{2}(n,k) = \Dimm{2}{n,k}\;.\]
\kmqed

\subsection{Product formula}
In $n$-torus, we know that \(\ds  \myZ{a_1,\ldots,a_n} = \san{1}{a_1} \cdots \san{n}{a_n}\) with 
$\ds \sum_{i=1}^{n} |a_i| = k$ are a basis of \(\ds P_{k}\).  
First, we observe the product formula when 1-torus. 
For each non-negative integers $m,n$ the followings hold:
\begin{align}
\noalign{Since \(\ds \sin m x \cos n x = \frac{1}{2} ( \sin(m+n)x + \sin
(m-n) x\)}
\sin m x \cos n x &= \frac{1}{2} \sin (m+n) x 
\;,\;\text{i.e., } \myZ{m}\myZ{-n} = \frac{1}{2}\myZ{m+n}
\quad\textbf{in}\quad P_{m+n}
\\
\noalign{\(\ds \sin m x \sin n x = \frac{-1}{2} ( \cos(m+n)x - \cos 
(m-n) x\) implies }
\sin m x \sin n x & = \frac{-1}{2} \cos(m+n)x 
\;,\;\text{i.e., } \myZ{m}\myZ{n} = \frac{-1}{2}\myZ{-(m+n)}
\quad\textbf{in}\quad P_{m+n}
\\
\noalign{and since \(\ds \cos m x \cos n x = \frac{1}{2} ( \cos(m+n)x + \cos 
(m-n) x\) implies } 
\cos m x \cos n x &= \frac{1}{2}  \cos(m+n)x 
\;,\;\text{i.e., }\myZ{-m}\myZ{-n} = \frac{1}{2}\myZ{-(m+n)}
\quad\textbf{in}\quad P_{m+n}
\end{align}
In general $n$-torus, the multiplication formula is as follows:
\[ \myZ{A} \myZ{B} = 
\prod_{i=1}^{n} \phi(A_i, B_i) \; 
\myZ{\psi(A_1, B_1), \ldots,  
\psi(A_n, B_n) } \]
where 
\(\phi, \psi\) are defined on 1-torus by 
\begin{alignat}{3}
\phi(a,b)& = 1, & \qquad \psi(a,b) &= a +  b  & \text{if } &\quad   a b = 0 
\label{eqn:our:ruleZ}
\\
\; 
\phi(a,b)&= -\sgn(a)/2, & \quad \psi(a,b) &= -(|a|+|b|) \qquad & \text{if
} &\quad  a  b > 0
\label{eqn:our:ruleP}
\\
\quad
\phi(a,b) &= 1/2 , & \quad \psi (a,b) &= |a|+|b|  &
 \text{if } &\quad a b < 0
\label{eqn:our:ruleN}
\end{alignat}

\kmcomment{
\begin{equation} 
(\lambda_i, C_i ) = \begin{cases}
\quad 
1,\qquad  A_i +  B_i & \text{if } A_i B_i = 0  \\
\; 
   -\sgn(A_i)/2,\quad -(|A_i|+|B_i|)  & \text{if } A_i  B_i > 0
\\
\quad
1/2 ,\quad |A_i|+|B_i|  
& \text{if } A_i B_i < 0
\end{cases}
\label{eqn:our:rule}
\end{equation}
}

\subsection{Differentiation}
The formula for differentiation is simple:
\begin{equation}
\frac{\pdel}{\pdel x_i } \myZ{a_1, \ldots, a_i, \ldots, a_n} = 
a_{i} \myZ{a_1, \ldots, - a_i, \ldots, a_n} \end{equation} 

\section{Proofs by another basis}
\subsection{Proof of Theorem \ref{thm:first}}  
Now, our bracket is given by 
\begin{align}
 \Pkt{\myZ{a_1,a_2}}{\myZ{b_1,b_2}}  \notag 
= &
\pdel_1 \myZ{a_1,a_2} \pdel_2 \myZ{b_1,b_2} -
\pdel_2 \myZ{a_1,a_2} \pdel_1 \myZ{b_1,b_2}\notag \\
=& 
a_1 b_2 \myZ{- a_1,a_2} \myZ{b_1,- b_2} -
a_2 b_1 \myZ{a_1,- a_2} \myZ{- b_1,b_2} \notag
\\
=& 
a_1 b_2 \phi(- a_1,b_1) \phi( a_2,- b_2) 
\myZ{\psi(- a_1,b_1),\psi( a_2,- b_2)} 
\label{eqn:our:bkt}
\\&\hspace{20mm}
-
a_2 b_1 
\phi( a_1, - b_1 ) \phi(- a_2, b_2) 
\myZ{\psi( a_1, - b_1 ),\psi(- a_2, b_2)} 
\notag 
\end{align}
In order to apply 
 \eqref{eqn:our:ruleZ}, 
 \eqref{eqn:our:ruleP} and  
 \eqref{eqn:our:ruleN} for the bracket \eqref{eqn:our:bkt}, we divide our discussion in
 several cases depending on the signature of $a_1 b_1$ or $a_2 b_2$ as
 below:
\begin{center}
\begin{tabular}{c|*{3}{c}} 
$\ds \null_{ \text{sign of }a_1 b_1 } \diagdown { }^ {\text{sign of }a_2
b_2 } $ &  $0$ & positive & negative \\\hline
$0$ & Case[1] & Case[2] & Case[3] \\
positive &         & Case[4] & Case[5] \\
negative &         &         & Case[6] \\
\end{tabular}
\end{center}

In the first three cases, we assume \(\ds a_1 b_1 =0\). 
If $a_1=b_1=0$ then our bracket is 0, so we may assume $\ds a_1 \ne 0$ 
and $\ds b_1 = 0$.  
Then our bracket is given by 
\begin{align} \text{ourBkt} = & 
a_1 b_2 \phi(- a_1,0) \phi( a_2,- b_2) 
\myZ{\psi(- a_1,0),\psi( a_2,- b_2)} 
= 
a_1 b_2  \phi( a_2,- b_2) 
\myZ{- a_1,\psi( a_2,- b_2)} 
\\
\noalign{ 
Case[1] \(a_2  b_2 = 0\): If \(b_2=0\) then ourBkt is 0, so we may
assume \(\ds a_2 =0, b_2 \ne 0\) and have 
}
\text{ourBkt} = &
a_1 b_2  \phi(0, - b_2) \myZ{- a_1,\psi( 0,- b_2)} 
= a_1 b_2  \myZ{- a_1,- b_2} \label{eqn:bkt:one}
\\
\noalign{ 
Case[2] \(a_2  b_2 > 0\): 
}
\text{ourBkt} =& 
a_1 b_2  \phi( a_2,- b_2) \myZ{- a_1,\psi( a_2,- b_2)} 
= 
a_1 b_2  \frac{1}{2} \myZ{- a_1, | a_2| + | b_2|} 
\\
\noalign{
Case[3] \(a_2  b_2 < 0\): 
}
\text{ourBkt} =& 
a_1 b_2  \phi( a_2,- b_2) \myZ{- a_1,\psi( a_2,- b_2)} 
= 
a_1 b_2  \frac{-\sgn{a_2}}{2} \myZ{- a_1, -| a_2| - | b_2|} 
\end{align}

We assume 
\(a_1  b_1 > 0\) in Case[4] or Case[5]. Then we have  
\begin{align}
\text{ourBkt} = &  
a_1 b_2 \frac{1}{2}  \phi( a_2,- b_2) 
\myZ{ |a_1|+|b_1|,\psi( a_2,- b_2)} 
-
a_2 b_1 
\frac{1}{2} \phi(- a_2, b_2) 
\myZ{ |a_1|+|b_1| ,\psi(- a_2, b_2)} 
\\
\noalign{ 
Case[4] \( a_2 b_2 > 0\): 
}
\text{ourBkt} = &  
a_1 b_2 (\frac{1}{2})^2   
\myZ{ |a_1|+|b_1|, |a_2|+|b_2|}
-
a_2 b_1 
(\frac{1}{2})^2 
\myZ{ |a_1|+|b_1| , |a_2|+|b_2| }
\notag
\\
=& 
\frac{1}{4} 
( a_1 b_2 -a_2 b_1 ) 
\myZ{ |a_1|+|b_1|, |a_2|+|b_2|}
\\
\noalign{ 
Case[5] \( a_2 b_2 < 0\): 
}
\text{ourBkt} = &  
a_1 b_2 (\frac{1}{2})^2 (-\sgn(a_2) )   
\myZ{ |a_1|+|b_1|, -|a_2|-|b_2|}
\notag \\& \hspace{20mm}
-
a_2 b_1 
(\frac{1}{2})^2 \sgn(a_2)    
\myZ{ |a_1|+|b_1| , -|a_2|+|b_2| }
\notag
\\
=& 
- \frac{1}{4} \sgn(a_2) 
( a_1 b_2 +a_2 b_1 ) 
\myZ{ |a_1|+|b_1|, -|a_2|-|b_2|}
\end{align}

In the last 
Case[6], 
 \( a_1 b_1 < 0\) implies  
\begin{align}
\text{ourBkt}
= &
a_1 b_2 (\sgn(a_1)/2)  \phi( a_2,- b_2) 
\myZ{- |a_1|-|b_1|,\psi( a_2,- b_2)} 
\\ & \hspace{20mm}
-
a_2 b_1 
(-\sgn( a_1)/2) \phi(- a_2, b_2) 
\myZ{ - | a_1| - |b_1 |,\psi(- a_2, b_2)} 
\notag 
\\
\noalign{ 
\( a_2 b_2 < 0\) yields 
}
\text{ourBkt}
= &
a_1 b_2 (\sgn(a_1)/2)  (- \sgn(a_2)/2) 
\myZ{- |a_1|-|b_1|, - |a_2|-|b_2|} 
\notag
\\ & \hspace{20mm}
-
a_2 b_1 
(-\sgn( a_1)/2) 
(\sgn( a_2)/2) 
\myZ{ - | a_1| - |b_1 |, - | a_2| - |b_2 |}
\notag 
\\
= &  - \frac{1}{4} \sgn(a_1)\sgn(a_2) (a_1 b_2 - a_2 b_1) 
\myZ{ - | a_1| - |b_1 |, - | a_2| - |b_2 |}
\end{align}

We conclude that on 2-torus we have \(\ds \myZ{a,b} = \frac{1}{a b}
\Pkt{ \myZ{-a,0}}{\myZ{0,-b}} \) with $a b \ne 0 $ from
\eqref{eqn:bkt:one}.  We call $\myZ{a,b}$ has a property (G) if $a b
\ne 0$. The cases above show that \(\ds \myZ{a,0} \) nor \(\ds \myZ{0,a}
\) does not obtain from the brackets \(\ds \Pkt{ \myZ{A} }{\myZ{B}} \). 

Thus, on \(\mathbb{T}^2\) the corank of the boundary operator  \(\ds
C_{ 1,w} = P_{w}
\mathop{\leftarrow}^{\pdel} 
C _ {2,w} = 
\sum_{i<j}^{i+j=w} P_{i}\wedge P_{j} (+ \Lambda^2 P_{w/2} \ 
\text{if}\  w\ \text{even})\) is 4 because the cokernel has a basis  \(    
\ds \myZ{\pm w,0} \) and  \( 
\ds \myZ{0,\pm w} \) when $w>0$.  
\kmcomment{
\ref{thm:first}
\begin{theorem}\label{B:thm:first}
On the symplectic 2-torus \(\ds \mT^2\), we consider the formal
Hamiltonian vector fields, we may identify with periodic ``polynomials''. 
We saturate the space with respect to the weight and have $m$-th chain
space $\ds \CS{w}{m}$ with the given weight $w$.  
We have a chain complex 
\[  0 \leftarrow \CS{w}{1} 
 \mathop{\leftarrow}^{\pdel} \CS{w}{2} 
 \mathop{\leftarrow}^{\pdel} \CS{w}{3}  
 \mathop{\leftarrow}^{\pdel} \cdots \]
The first Betti number of this chain complex is 4 for each positive
weight $w$. 
\end{theorem}
}

\subsection{Proof of Theorem \ref{thm:dim4}} 
\kmcomment{ 
\ref{thm:dim4}
\begin{theorem}
The second difference of the sequence of 
the first Betti number 
of the weight $w$ on $\ds \mT^{4}$ 
is an arithmetic progression with common difference 16.  
More precisely, 
the first Betti number for the weight $w>0$ is given by  
\( 4^2 w -8 \).
\end{theorem}
\textbf{Proof:} 
}
Since 
\(\ds  \myZ{\ds a_1,a_2,a_3,a_4} = 
 \myZ{a_1,a_2,0,0}\myZ{0,0,a_3,a_4}\), we may denote   
\(\ds \myZ{A,A'} = \myZ{A} \myZ{A'} \) with no confusion.    
Since 
the Poisson bracket satisfies the Leibniz rule, 
we have \kmcomment{using the notation \(\ds [1,1] = \Delta\),} 
\begin{align}
 \Pkt{\myMM{A,A'}{P,P'}}{\myMM{B,B'}{Q,Q'}} \label{B:dim4:mk:master}
=  & 
 \Pkt{\myMM{A}{P}}{\myMM{B}{Q}}^{(2)} \myMM{A'}{P'}\myMM{B'}{Q'}
 + \myMM{A}{P} \myMM{B}{Q} 
 \Pkt{\myMM{A'}{P'}}{\myMM{B'}{Q'}}^{(2)} 
\end{align} 
where \(\ds {\Pkt{\cdot}{\cdot}}^{(2)}\) means the Poisson bracket in
2-dimensional torus.  
In \eqref{B:dim4:mk:master}, put \(\ds B'=O\) or \(\ds B=O\).   
Then because of 
 \(\ds \myMM{B'}{Q'} \) or \(\ds \myMM{B}{Q} \)  is constant $1$, we have 
\begin{align}
 \Pkt{\myMM{A,A'}{P,P'}}{\myMM{B,O}{Q,O}} 
= & \Pkt{\myMM{A}{P}}{\myMM{B}{Q}}^{(2)} \myMM{A'}{P'} \label{B:cok:one}\\
\Pkt{\myMM{A,A'}{P,P'}}{\myMM{O,B'}{O,Q'}} 
= & 
 \myMM{A}{P} \Pkt{\myMM{A'}{P'}}{\myMM{B'}{Q'}}^{(2)} \label{B:cok:two}
\end{align} 

From 
\eqref{B:cok:one} and \eqref{B:cok:two}, we have a diagram
\[
\begin{array}{c||c|c} 
{ }_{ V_{i,2}} \diagdown^{V_{j,2}} & \text{out of image in }\mT^2 & \text{image in
}\mT^2 \\\hline\hline
\text{out of image in }\mT^2 & \text{out of image in }\mT^4  
& \text{image in }\mT^4 \\\hline
\text{image in }\mT^2 & \text{image in }\mT^4 & \text{image in } \mT^4
\end{array}
\]

Here we introduce the notation $\crk{w}{2n}$ the corank of the boundary
operator \(\pdel\) to \(\CS{w}{1}\) where \(\CS{w}{1}\) is the 1st chain
space with the weight $w$, that is,  the $w$-homogeneous "polynomial"
space \(P_{w}\) on \(\mT^{2n}\).  

We know that $\crk{0}{2}=1$ and $\crk{w}{2}=4$ for $w>0$ in Theorem
\ref{thm:first}. The table above implies that    
\begin{align*} 
 \crk{w}{4} 
=& 
 \crk{0}{2} \crk{w}{2} \\& 
 + \sum_{i=1}^{w-1} \crk{i}{2} \crk{w-i}{2} 
 + 
  \crk{w}{2} \crk{0}{2} 
  \\
=&  4 + 4^2 (w-1) + 4 = 4^2 w -8 
\\
\crk{0}{4} =& 1 
\end{align*}
\kmqed

\subsection{Proof of Theorem \ref{thm:gen}}

Observing the discussion in \(\ds \mT^4\) case, we show the same
discussion works well for general even dimensional torus \(\ds
\mT^{2n}\) case, too.  
\kmcomment{
 \ref{thm:gen}
\begin{theorem}
Let \(\ds \crk{w}{2n}\) be the dimension of the first homology group of \(\ds
 V_{w,2n} \mathop{\leftarrow}^{\pdel } \CS{w}{2}\). Then we have the
 next recursive  formula:  
\begin{align*} 
 \crk{w}{2n} 
=&  4 + 4 \sum_{i=1}^{w}  \crk{i}{2n-2} - 3 \crk{w}{2n-2} 
\\
\crk{0}{2n} =& 1 
\end{align*} 
\end{theorem}
\textbf{Proof:} 
}
According to the symplectic structure, we separate $2n$-sequence to $n$
pairs and  
\(\ds \myMM{\myUP{A}{1}\ldots \myUP{A}{n}}
              {\myUP{P}{1}\ldots \myUP{P}{n}}\) are our basis.  
To avoid complicated expression, we abbreviate the first $(n-1)$
sequences 
\(\ds \myUP{A}{1}\ldots \myUP{A}{n-1}\) by \(A\).  
Then the Poisson bracket satisfies 
\begin{subequations}
\begin{align} 
\Pkt{ 
 \myMM{A\;\myUP{A}{n}}
              {P\;\myUP{P}{n}}}
{ \myMM{B\; \myUP{B}{n}}
              {Q\;\myUP{Q}{n}}}_{(2n)}    
=& 
\Pkt{ 
 \myMM{A} {P}}
{ \myMM{B }{Q} } _{(2n-2)} 
 \myMM{\myUP{A}{n} }{\myUP{P}{n}} 
 \myMM{\myUP{B}{n} }{\myUP{Q}{n}} 
\label{B:shiki:1}
              + 
 \myMM{A} {P} 
 \myMM{B} {Q} 
 \Pkt{ \myMM{\myUP{A}{n}}{\myUP{P}{n}}}
{ \myMM{\myUP{B}{n}}{\myUP{Q}{n}}}_{(2)}
 \end{align}
\end{subequations}
\(\ds \myUP{B}{n}=O\) implies 
\( 
 \myMM{\myUP{B}{n}}{\myUP{Q}{n}} = 1\) and 
\begin{subequations}
\begin{align} 
\Pkt{ 
 \myMM{A\;  \myUP{A}{n}}
              {P\; \myUP{P}{n}}} 
{ \myMM{B\; O} {Q\; O}} 
=& 
\Pkt{ 
 \myMM{A } {P }} 
{ \myMM{B } {Q }} 
 \myMM{\myUP{A}{n} }
              {\myUP{P}{n} } 
\label{B:shiki:2:2}
\end{align}
\end{subequations}
Now multiplying 
\(\ds\myZ{ \myUP{B}{n}}\) to  
\(\ds \myZ{ A\; \myUP{A}{n}} \), 
we have 
\begin{subequations}
\begin{align} 
&
\Pkt{
\myMM{A\; \myUP{A}{n} } 
          {P\; \myUP{P}{n} }
\myMM{\myUP{B}{n} } 
          {\myUP{Q}{n} }
          } 
    { \myMM{B \;  O} {Q\; O}}  
= 
\Pkt{ 
  \myMM{A } {P }} 
{ \myMM{B } {Q }}
 \myMM{\myUP{A}{n} } {\myUP{P}{n} } 
 \myMM{\myUP{B}{n} } {\myUP{Q}{n} } 
\label{B:shiki:3}
 \end{align}
\end{subequations}
Subtracting \eqref{B:shiki:3} from \eqref{B:shiki:1}, we have 
\begin{subequations}
\begin{align}
&
\Pkt{ 
 \myMM{A\; \myUP{A}{n}}
              {P\;  \myUP{P}{n}}} 
{ \myMM{B\; \myUP{B}{n}}
              {Q\; \myUP{Q}{n}}} -  
\Pkt{ 
\myMM{A\;  \myUP{A}{n}} {P \; \myUP{P}{n} }
\myMM{ \myUP{B}{n} } {\myUP{Q}{n} }
} 
{ \myMM{B\; O} {Q\; O}} 
= 
 \myMM{A } {P } 
 \myMM{B } {Q } 
\Pkt{ \myMM{\myUP{A}{n}}{\myUP{P}{n}}}
{ \myMM{\myUP{B}{n}}{\myUP{Q}{n}}} \label{B:shiki:4:2} 
 \end{align}
\end{subequations}

\eqref{B:shiki:2:2} and  
\eqref{B:shiki:4:2} imply the diagram
\[
\begin{array}{c||c|c} 
V_{i,2n-2} \backslash V_{j,2} & \text{out of image in }\mT^2 & \text{image in
}\mT^2 \\\hline\hline
\text{out of image in }\mT^{2n-2} &  
\text{out of image in }\mT^{2n}   & \text{image in }\mT^{2n} \\\hline
\text{image in }\mT^{2n-2} & \text{image in }\mT^{2n} & \text{image in }
\mT^{2n}
\end{array}
\]
and we see that for $w>0$ 
\begin{align*} 
 \crk{w}{2n} 
=& 
 \crk{0}{2n-2} \crk{w}{2} \\& 
 + \sum_{i=1}^{w-1} \crk{i}{2n-2} \crk{w-i}{2} 
 \\& 
 + \crk{w}{2n-2} \crk{0}{2} 
  \\
=&  4 + 4 \sum_{i=1}^{w-1}  \crk{i}{2n-2} + \crk{w}{2n-2} 
\\
\crk{0}{2n} =& 1\;.  
\end{align*}
This theorem contains Theorem \ref{thm:dim4} when  
$2n=4$.
\kmqed

\kmcomment{
\ref{cor:dim6}
\begin{kmCor}
On \(\ds \mT^6\), 
 \(\ds \crk{0}{6}  = 1\) and for \(\ds w>0\) we have   
\begin{equation}  
 \crk{w}{6}  = 4 ( 8 w^2 - 12 w + 7)\; .  
 \end{equation} 

On \(\ds \mT^{2n+2}\), it holds 
\begin{equation}\ds 
\sum_{k=0}^ n (-1)^{k} \binom{n}{k} 
 \crk{w-k}{2n+2}  = 2^{2n+2} \label{eqn:recursive2} \end{equation}  
for \( w > n\). 
\end{kmCor}

\textbf{Proof of Cor \ref{cor:dim6} :} 
Combining Theorem \ref{thm:dim4} and Theorem \ref{thm:gen}, direct
computation yields the formula 
\eqref{eqn:dim6}.  
We have   
\begin{align*}
& 
\sum_{k=0}^{1} (-1)^{k} \binom{1}{k} \crk{w-k}{4}  =  
\sum_{k=0}^{1} (-1)^{k} \binom{1}{k} ( 16(w-k) - 8) \\ =&   
\sum_{k=0}^{1} (-1)^{k} \binom{1}{k} ( 16w- 8)  + 
\sum_{k=0}^{1} (-1)^{k} \binom{1}{k} ( - 16k )    =
0 + 16 = 2^{4} 
\end{align*}
Thus, \eqref{eqn:recursive2} holds for $2n=2$.  
\kmcomment{
Also,  
\begin{align*}
& 
\sum_{k=0}^{2} (-1)^{k} \binom{2}{k} 
 \crk{w-k}{6}  \\=&  
 \crk{w}{6}  - 2   
 \crk{w-1}{6}  +   
 \crk{w-2}{6}   = 2^6\;.  
\end{align*}
Thus, \eqref{eqn:recursive2} holds for $n=2$.  
}
Assume \eqref{eqn:recursive2} holds on \(\ds \mT^{2n+2}\) for general $n$.   
First we confirm that 
\begin{align} & 
\sum_{k=0}^{n+1} (-1)^k \binom{n+1}{k} \crk{w-k}{2n+2}
\notag
\\ =& 
\sum_{k=0}^{n+1} (-1)^k ( \binom{n}{k} 
+\binom{n}{k-1})  
\crk{w-k}{2n+2} = 2^{2n+2}-2^{2n+2} = 0\;.  
\label{eqn:hojo2}
\end{align}
Now we prove 
\eqref{eqn:recursive2} holds on \(\ds \mT^{2n+4}\).  
\begin{align*} & 
\sum_{k=0}^{n+1} (-1)^k \binom{n+1}{k} \crk{w-k}{2n+4} \\
=& \sum_{k=0}^{n+1} (-1)^k \binom{n+1}{k} \left( 4+ 4 \sum_{i=1}^ {w-k}
\crk{i}{2n+2} - 3 \crk{w-k}{2n+2} \right) 
\\ 
=& 4 \sum_{k=0}^{n+1} (-1)^k \binom{n+1}{k}  - 3 
 \sum_{k=0}^{n+1} (-1)^k \binom{n+1}{k} 
\crk{w-k}{2n+2} \\& 
+ 
 4 \sum_{k=0}^{n+1} (-1)^k \binom{n+1}{k}  \sum_{i=1}^ {w-k}
\crk{i}{2n+2} \\
=& 
 4 \sum_{k=0}^{n+1} (-1)^k \binom{n+1}{k}  \sum_{i=1}^ {w-k}
\crk{i}{2n+2} \\
\noalign{because of general fact about alternating sum, and 
\eqref{eqn:hojo2}. Changing the order of summation, we have }
=& 4 \sum_{i=1}^{w-n-1} \sum_{k=0}^{n+1} (-1)^k
\binom{n+1}{k}\crk{i}{2n+2} 
+ 4 \sum_{i=w-n}^{w} \sum_{k=0}^{w-i} (-1)^k \binom{n+1}{k} \crk{i}{2n+2} 
\\
=& 0  
+ 4 \sum_{i=w-n}^{w} \sum_{k=0}^{w-i} (-1)^k \binom{n+1}{k} \crk{i}{2n+2} 
= 
 4 \sum_{i=w-n}^{w} \sum_{k=0}^{w-i} (-1)^k ( \binom{n}{k} 
  +\binom{n}{k-1} ) \crk{i}{2n+2} 
  \\
= & 
 4 \sum_{i=w-n}^{w} (-1)^{w-i} \binom{n}{w-i}\crk{i}{2n+2} = 4 \cdot
 2^{2n+2}
 = 2^{2n+4} \;. 
\end{align*}
\kmqed
}

\end{appendices}

\end{document}